%% file: LevelSetCancerInvasion.tex
\documentclass[12pt]{article}
\usepackage[square,sort,comma,numbers]{natbib}
\usepackage[english]{babel}
\usepackage[utf8x]{inputenc}
\usepackage{amsmath,amsfonts}
\usepackage{color}
\usepackage{subfig}
\usepackage{float}
\usepackage{hyperref}
\usepackage{graphicx}
\usepackage{float}
\usepackage[noend]{algorithmic}
\usepackage[ruled]{algorithm}

\newcommand{\dy}{\rm{d}y}
\newcommand{\dz}{\rm{d}z}
\newcommand{\dt}{\rm{d}t}
\newcommand{\dvar}[1]{\frac{\partial #1}{\partial t}}
\newcommand{\dvarp}[1]{\frac{\partial #1}{\partial \tau}}
\newcommand{\mh}{\widehat m}

\newcommand{\mhhn}{\widehat m_{h,n}}
\newcommand{\T}{\mathcal M_h}
\newcommand{\VMh}{U_M^h}
\newcommand{\VTh}{U_T^h}
\newcommand{\Vmh}{U_m^h}
\newcommand{\eps}{\varepsilon}
\newcommand{\cv}{c_{\rm vel}}

\newtheorem{remark}{Remark}[section]
\newtheorem{assumption}{Assumption}[section]

\usepackage{amssymb}
\newcommand{\tho}[1]{\textcolor{black}{#1}}

\newcommand{\dtr}[1]{\textcolor{black}{#1}}
\usepackage[affil-it]{authblk}

\title{A level-set approach for a multi-scale cancer invasion model}

\author{Thomas Carraro$^{1}$\thanks{thomas.carraro@iwr.uni-heidelberg.de}, Sven E. Wetterauer$^{1}$,\\ Ana Victoria Ponce Bobadilla$^{1}$, Dumitru Trucu$^{2}$}
\affil{$^1$Institute for Applied Mathematics, Heidelberg University,\\ 69120 Heidelberg, Germany}
\affil{$^1$Interdisciplinary Center for Scientific Computing (IWR),\\ Heidelberg University, 69120 Heidelberg, Germany}
\affil{$^2$Division of Mathematics, University of Dundee,\\ Dundee, DD1 4HN, United Kingdom}

\begin{document}

\maketitle

%\begin{frontmatter}



%\cortext[cor1]{Corresponding author.}
\begin{abstract}
\dtr{Central to the quest for a deeper understanding of the cancer growth and spread process, the naturally multiscale character of cancer invasion demands appropriate multiscale modelling and analysis approach. The cross-talk between the tissue scale (macro-scale) cancer cell population dynamics and the cell-scale (micro-scale) proteolytic molecular processes along the tumour boundary plays a particularly important role within the invasion processes, leading to dramatic changes in tumour morphology and influencing the overall pattern of cancer spread. }

\dtr{Building on the multiscale moving boundary framework proposed in Trucu et al. \cite{trucu2013multiscale} (Multiscale Model. Simul 11(1): 309-335), in this work we propose a new formulation of this process involving a novel derivation of the macro-scale boundary movement law based on micro-dynamics, involving a transport equation combined with the level-set method. This is explored numerically in a novel finite element macro-micro framework based on cut-cells.}
\end{abstract}

\section{Introduction}

Involving a wide range of cross-related processes occurring on several spatio-temporal scales, cancer cell invasion of the human is one of the hallmarks of cancer \cite{hanahan2011hallmarks}, playing a crucial role in the overall development a spread of a growing malignant tumour. Taking advantage of the heterotypic character of the tumour microenvironment (which includes immuno-inflamatory cells, stromal cell, fibroblasts, endothelial cells, macrophages), complex molecular processes facilitate intense interactions between the cancer cell population and the extracellular matrix (ECM)\cite{hanahan2011hallmarks,Egeblad_et_al_2010,qian2010macrophage,Joyce_Pollard_2009,Kalluri_Zeisberg_2006}. These \tho{interactions} lead to a cascade of specific developmental patterns and behaviours of the growing tumours, most notable stages including the degradation of the ECM, the local progression of the tumour, followed by the tumour angiogenesis process and the subsequent metastatic spread of the cancer cells in the human body. 

Focusing on local tumour progression, the alteration and remodelling of the ECM by the matrix degraded enzymes (MDEs) such as matrix metalloproteinases MMPs or the urokinase plasminogen activator (uPA) play a key role. Alongside cell-adhesion and multiple taxis processes (including haptotaxis and chemotaxis), the matrix degrading enzymes processes degrade various components of the surrounding ECM that leads to further tumour progression. However, as the full mechanisms involved in these complex processes is yet to be deciphered biologically, over the past two decades or so cancer invasion received extensive mathematical modelling attention, in which systems of reaction-diffusion-taxis partial differential equations \cite{Adam_1986,Anderson_et_al_2000,Byrne_Chaplain_1996,Chaplain_et_al_2006,Gatenby_Gawlinski_1996,Greenspan_1976,Perumpanani_et_al_1996,Perumpanani_et_al_1998,Webb_et_al_1999,Andasari_et_al_2011,Anderson_2005,Byrne_et_al_2001,Chaplain_Lolas_2005,Deakin_Chaplain_2013} as well as nonlocal integro-differential systems \cite{Armstrong_et_al_2006,Chaplain_et_al_2011,Domschke_et_al_2014,Gerisch_Chaplain_2008} were derived and proposed to deepen the understanding, validate and create new experimental hypothesis. 
Furthermore, to capture various heterotypic aspects and related processes within tumour invasion, several multiphase models based on the theory of mixtures \cite{Byrne_Preziosi_2003,Chaplain_et_al_2006b,Frieboes_et_al_2010,Preziosi_Tosin_2009,Wise_et_al_2008,Wise_et_al_2011} were derived (by exploring the mass and momentum balances as well as the inner multiphase constitutive laws). 

A particularly important role in cancer invasion is played by the MDEs (such as the MMPs) that are secreted from the outer proliferating rim and released within the tumour peritumoural microenvironment. This gives the cancer invasion a moving boundary character, and to that end several level-set approaches were recently proposed to study the tumour progression both in  homogeneous environments \cite{Frieboes_et_al_2006,macklin2005evolving,macklin2006improved,macklin2007nonlinear,zheng2005nonlinear} and in complex heterogeneous tissues \cite{macklin08}. 

Despite recent advances, the multiscale modelling of the processes involved in cancer invasion remains an open problem. Although this is a truly multiscale process, most mathematical models were offering a one-scale perspective, whether that is from a purely macro-scale (tissue scale) or an exclusively micro-scale (cell-scale) stand point. 
However, recently a novel 2D multiscale moving boundary modelling platform for cancer invasion was proposed in \cite{trucu2013multiscale}. This explores in an integrated manner the tissue-scale cell population dynamics and relevant cell-scale molecular mechanics together with the permanent link between these two biological scales. This addresses directly the dynamics of the MDEs proteolytic processes occurring at the tumour boundary (i.e., at the invasive edge of the tumour) that are sourced from within outer proliferating rim of the tumour and facilitate the complex molecular transport and ECM degradation within the peritumoural region.  The tissue-scale progression of tumour morphology is captured here in a multiscale moving boundary approach where the contribution arriving from the cell-scale activity to the cancer invasion pattern is realised by the micro-scale MDEs dynamics (occurring along the tumour invasive edge), which, for its part, is induced by the cancer macro-dynamics. This was recently applied to the extended context in which, rather than the MMPs dynamics, the uPA is considered as the proteolytic system, and has led to biologically relevant results \cite{Peng2017}.

In this work we present a new formulation of the link between the two scales presented in \cite{trucu2013multiscale}. The new model is based on a level set approach in which the moving domain is defined as the zero level of a level set function. The reason for this choice of the problem setting is twofold: on one hand, all components of the model can be described by partial differential equations at the continuum level allowing the complete separation between modelling and discretization; on the other hand, it is better suited for an extension to the three-dimensional case since the formulation of all components of the model is dimension independent and the use of a dimension independent implementation of the discretization, like in our case using the finite element method (FEM) package deal.II \cite{dealii:2017},  facilitates the realization of the code.

The level set method was first introduced in \cite{osher1988fronts} for tracking moving interface with complex deformations. 
This method was developed starting from the notion of weak solutions for evolving interfaces.
The main aspect of this method is the fact that an interface or a domain boundary is defined through the embedding of the interface as the zero level set of a higher dimensional function. Furthermore, the velocity of the interface is also embedded to the higher dimensional function. \tho{We avoid handling a sharp interface, i.e.\ a lower dimensional manifold in the computational domain, but the velocity needs to be extended from the interface to the rest of the domain}. While the original setting, with the sharp interface, poses several numerical difficulties due to its Lagrangian approach, the later setting, using an Eulerian approach, can exploit techniques developed for hyperbolic problems.

Other works have presented a level set approach for moving the tumour interface. In \cite{macklin2005evolving} the authors use a level set function to define the boundary of a tumour mass and extend the velocity orthogonally to the interface using a filter technique to damp numerical noise coming from the extension procedure. The velocity at the interface is defined as a function of the gradient of a computed quantity (the pressure). This work nevertheless does not link different model scales. 
In \cite{zheng2005nonlinear} an adaptive finite element combined with a level set approach is used to solve a model that considers tumour necrosis, neo-vascularization and tissue invasion. The model is composed of a continuum part and a hybrid continuum-discrete part. The velocity of the interface is the cell velocity. Therefore, \tho{the velocity does not need to be extended into} the neighbourhood of the interface.
In \cite{macklin2009multiscale} a level set approach with a ghost-cell method is applied to tumour growth of glioblastioma. The velocity of the interface depends on solutions of linear and nonlinear equations with curvature-dependent boundary conditions. Since the velocity is only defined at the interface the authors extend it beyond the interface and use a narrow band/local level technique to update the interface velocity and level set function only in the vicinity of the interface. Our approach uses a level set method with an extension of the velocity. While the coupling between the macroscopic and microscopic scales was originally introduced in \cite{trucu2013multiscale} considering a Lagrangian approach to move the nodes of the discrete approximation of the interface, we introduce here a continuous link between the two scales defined by the velocity at the interface at the continuum level. This changes the formulation of the multiscale coupling that goes with the definition of the velocity starting from heuristic arguments. The scope of this work is to present as a whole the new formulation of the tumour invasion model explaining the possible advantages that this approach can have for future developments and the numerical aspects that need further attention and further development.

The paper is organized as follows. In Section \ref{problem description} we state the problem setting and describe the different components of the model: the macroscopic and microscopic components and the description of the moving boundary. In Section \ref{weak form} we introduce the weak formulation of the model which is needed for the approximation of the continuum problem with a finite element method. We introduce the discretization of the problem using cut-cells for the approximation of the cancer region domain.
In Section \ref{num results} we present some numerical results showing the interplay of the different parts of the multiscale model. Finally we present an outlook and some concluding remarks in Section \ref{conclusions}.
\section{Problem description}
\label{problem description}
We present a \dtr{two-scale} model for cancer invasion that \dtr{links through a double feedback loop the dynamics occurring at two different spatial scales explored by the following two modelling components, namely}: a macroscopic \dtr{component} describing the population of cancer cells and extracellular matrix \dtr{at tissue-scale} and a microscopic \dtr{component} describing the dynamics of a \dtr{generic matrix-degrading} enzyme \dtr{molecular} population \dtr{at cell-scale}. Both scales are considered at the continuum level. We therefore assume that possible stochastic effects, in regions where the continuum assumption is not more valid, can be neglected. Nevertheless, our aim is to derive a flexible numerical framework that would allow to extend the model with a stochastic part (e.g.\ at the interface of the domain) leading to a hybrid formulation, if needed.

The cancer cells and extracellular matrix are modelled in a\dtr{n invading} domain $\Omega(t)$ that changes \dtr{its size and morphology} in time during the invasion process \dtr{within a reference maximal tissue cube $Y$}.  \dtr{Its boundary $\partial \Omega(t)$ will also be referred to as {\em interface}} because \dtr{this is regarded here as} an interface separating the region with zero cancer cells and a given distribution of ECM from a region with a distribution of cancer cells and ECM that satisfies the macroscopic equations. \dtr{Finally, in appropriate cell-scale neighbourhoods of the interface $\partial \Omega(t)$} points, a microscopic problem describing the \dtr{cross-interface} transport of \dtr{matrix-degrading} enzymes \dtr{is considered and accounted for in order to determine the law for macro-scale boundary movement.}
\begin{assumption}[Scale separation]
\label{ass:scale separation}
We assume scale separation in space between the two components of the model, i.e.\ the macroscopic populations of cancer cells and ECM and the microscopic population of \dtr{matrix degrading }enzyme molecules. 
\end{assumption}
\dtr{The} characteristic length $L$ for the macroscopic part of the model \dtr{relates} to the diameter of the cancer region and \dtr{is considered here as in \cite{Anderson_et_al_2000,Gerisch_Chaplain_2008}, ranging between $0.1$cm and $1.0$cm. Further,} the characteristic length $\ell$ of the microscopic part related to the region where \dtr{the matrix degrading enzymes are spatially transported is considered to be of the order of $10^{-3}$cm, \cite{Mumenthaler:2013}.}
The ratio between the scales is denoted $\eps=\ell/L$.
Therefore, the enzyme population's dynamics can be described in a \dtr{bundle of} microscopic domain\dtr{s} $\eps Y$ that \dtr{are} obtained by scaling a reference \dtr{tissue cubic} domain $Y$ by a factor $\eps$ \dtr{and whose union provide a cell-scale neighbourhood for the interfacial points in $\partial \Omega(t)$}.
Due to spatial scale separation\dtr{,} we consider a microscopic problem \dtr{at} each point $x \in \partial \Omega(t)$ \dtr{(i.e., at each point of the macroscopic interface) on the corresponding $\varepsilon Y$ micro-domain centred at $x$}. \dtr{Further assuming for convenience that the maximal cube $Y$ is centred at origin of the space,} the micro-scale coordinates $y$ of the microscale problem \dtr{on a $\varepsilon Y$ centred at $x_{0}$} are obtained by \dtr{appropriate} scaling \dtr{and translation of }the macroscopic coordinates\dtr{, namely as} $y = \dtr{x_{0}+}\eps (x\dtr{-x_{0}})$. 
This approach is similar to the heterogeneous multiscale method for interface dynamics presented in \cite{cheng2003heterogeneous}, but we do not consider an upscaling process.

The first step of this work is to describe the whole coupled problem at the continuum level and then in a second step to discretize it. In particular, we consider a semidiscretization in time by the implicit \dtr{E}uler method and the finite elements method for the discretization in space. The reasons for this choice will be clar\dtr{ified} later in the section dedicated to the discretization.

\begin{figure}[!h]
\begin{center}
	\includegraphics[width=.7\textwidth]{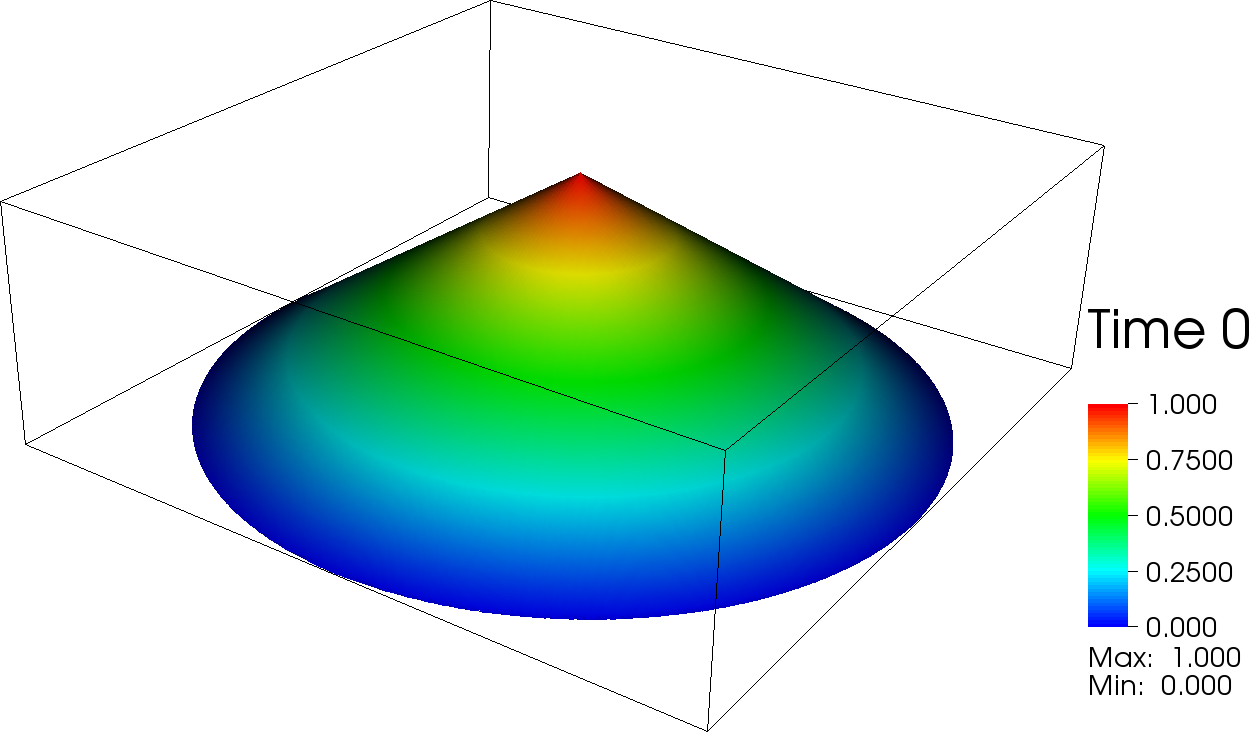}
\end{center}
\caption{Initial distribution of cancer cells}
\label{fig:init cancer}
\end{figure}

Schematically the two scales are coupled \dtr{by the following double feedback loop}:
\begin{itemize}
\item The \dtr{\emph{top-down}} macroscopic-to-microscopic coupling \dtr{for the microscopic problem at any boundary point $x\in \Omega(t)$} is done \dtr{via the source of matrix-degrading enzymes}  \dtr{which is induced by the macro-dynamics and is formed as a collective contribution of the cancer cells that arrive during the macro-dynamics within an appropriate distance from $x$}, see \eqref{F};
\item The \emph{bottom-up} microscopic-to-macroscopic coupling is done by defining the velocity of the macroscopic interface in dependence of the enzymes concentration in the microscopic domains, see \eqref{eq:velocity}.
\end{itemize}
Therefore, the rate of cancer cells invasion into the surrounding tissues is defined by the velocity of the interface that depends on the solution of the microscopic enzyme \tho{dynamics}.
The interface is described by the zero-level of a level-set function that has to be initialized with the initial distribution of the cancer cells, see Figure \ref{fig:init cancer} and formula \eqref{initial conditions}, and that is moved by a transport equation. 
%At the continuum level, it is sufficient to define the velocity only at the interface to fully describe the cancer invasion process. We introduce later an extension of the velocity in the vicinity of the interface needed in the weak formulation that we present later for the computations.

To account for all interactions between the different parts of the problem we introduce the domain $Y\subset \mathbb R^2$ that is assumed to be sufficiently large such that the complete dynamics happen inside it. This domain is used at the continuum level to define the region in which the transport equation is solved. Later it is used for the discretized problem to define a region where the finite elements become active if their intersection with the cancer region is not empty.

\dtr{In the following we proceed with} the \dtr{multiscale} model description in three parts: macroscopic, microscopic and transport component. \dtr{At macro-scale, t}he model considers cancer cells and extracellular matrix (ECM) interaction. 
\subsection{Macroscopic model component}\label{section21_macro_model}
Let $c(x,t)$ and $v(x,t)$ denote the cancer and the extracellular matrix distribution\dtr{s at $(x,t)\in \Omega(t)\times \tho{(0,T)}$,} respectively.  \dtr{Proceeding as in \cite{trucu2013multiscale}, the dynamics at} macroscopic scale \dtr{is given by} the following \dtr{PDE system}: 
\begin{align}
\label{eq:macrodynamics 1}
\frac{\partial c}{\partial t} &= \overbrace{ D_1\Delta c}^\text{Random motility}-\overbrace{\eta\nabla\cdot(c\nabla v)}^{Haptotaxis}+\overbrace{\mu_1(v)\, c\, (1-c-v)}^\text{Proliferation},  \\
\label{eq:macrodynamics 2}
    \frac{\partial v}{\partial t} &= -\underbrace{\alpha cv}_\text{Degradation}+\underbrace{\mu_2(1-c-v)}_\text{ECM Remodelling },
\end{align}
where $D_1$ is the diffusion coefficient for the cancer cells, $\eta$ is the advection coefficient, $\mu_1$ is the proliferation coefficient, $\alpha$ a degradation coefficient and $\mu_2$ a coefficient for the remodelling of ECM. The typical values of these coefficients \tho{are} reported in Table \ref{tab:num setting}. \dtr{Except for }$\mu_1$, \dtr{all other coefficients} are considered constant. \dtr{However, due to the direct dependency of the mitotic process on the presence of ECM \cite{chiarugi2008,pickup2014,tilghman2010}, the proliferation coefficient $\mu_1$ is a function of $v$ which decreases monotonically with the amount of ECM, exploring the full range of values between 1 and 0, as $v$ varies from its maximal levels, which are ideal for proliferation, to 0-level ECM regions, where cell-death and necrosis occurs due to lack of anchorage and nutrients.} 

It is assumed that the cancer cells are zero outside $\Omega(t)$ and that there is no flux of cells through the boundary $\partial \Omega(t)$. \dtr{Furthermore, under the presence of these boundary and initial conditions, for the case of constant proliferation rate $\mu_{1}$, the results in \cite{szymanska2009,winkler2010,hillen2013} explore the local and global existence of system \eqref{eq:macrodynamics 1}-\eqref{eq:macrodynamics 2}.}
The initial distribution of cancer cells $c_0(x)$ and extracellular matrix $v_0(x)$ are given in the larger domain $Y$.

\dtr{Finally, we note that while this macro-dynamics occurs at tissue-scale, the macroscale cancer invasion is completely controlled by the movement of the tumour interface which is reg\dtr{u}lated by the macroscopic/microscopic interaction that defines the interface velocity, whose details will be discussed in the next subsections.}

\dtr{At this stage, it is important to note that t}he macroscale part \eqref{eq:macrodynamics 1}-\eqref{eq:macrodynamics 2} of the model \dtr{can be regarded as} a singularly perturbed diffusion-reaction-transport problem of the type
\begin{align}
\frac{\partial c}{\partial t} - D_1 \Delta u + \eta \nabla \cdot(c \nabla v) - \mu_1(v) g_1(c,v) &= 0,\\
\frac{\partial v}{\partial t} - \gamma \Delta v + \alpha c v - \mu_2 g_2(c,v) &= 0
\end{align}
with $\gamma=0$. In particular, \dtr{we observe that} the second equation is reaction dominated and, in this context, it is well known that this kind of problems could present one or more boundary layers \cite{Surla:1997,John:2007}. The thickness of the boundary layers give these problems a multiscale character \dtr{in its own right}. 

Classical techniques to solve convection/reaction dominated problems are discontinuous and continuous interior penalty methods \cite{Burman:2005}.
Other techniques are streamline-upwind/Petrov-Galerkin formulations \cite{Fries:2004}.
Furthermore, to overcome numerical instabilities other methods have been developed like the Galerkin enriched finite element method (GEM) \cite{Fernando:2012}.
A typical problem in case of internal or boundary layers is the missing a priori knowledge of the position of such layers. Nevertheless, using a range of model parameters limited to the values used in literature for our model, we have observed that the solution component $v$ has a mild boundary layer, due to the degradation of the ECM.
Therefore, for the simulations shown here no stabilization technique was necessary and the mesh has been refined for accuracy purposes and not for stability issues.

In the next section we introduce the part of the model used to update the domain in time.

\subsection{Macroscopic time-dependent domain}
\dtr{As mentioned already above, the movement of the tumour interface is directly governed by the matrix degrading enzymes (MDE) dynamics occurring in a cell-scale neighbourhood of the tumour interface $\partial \Omega(t)$}. \dtr{The pattern of degradation of the peritumoural ECM by the advancing front of MDEs drives the invasion of the tumour cells in the surrounding tissues and determines the movement of the tumour boundary $\partial \Omega(t)$. Therefore, } the movement of the time-dependent macro domain $\Omega(t)$ is \dtr{enabled} by a velocity field defined on the points of the interface $x\in \partial \Omega(t)$, \dtr{which is determined by the micro-dynamics occurring on a small micro-domain $\varepsilon Y$ centred at $x$}.
\dtr{Hence,} the velocity field \dtr{that arises this way on $\partial\Omega(t)$} depends on the micro\dtr{-dynamics MDE molecular distribution $m(y,\tau)$ over an appropriate micro-spatio-temporal domain $\varepsilon Y\times (0,\Delta T)$ (which will be detailed in Section} \ref{sec:micro}). \dtr{Therefore, we will denote this velocity field by $V(m)$ and remark at this stage that this establishes a bottom-up feedback link from micro- to macro-dynamics, being directly responsible for the movement of the tissue scale tumour boundary $\partial \Omega(t)$}.

Since \dtr{$V(m)$} is defined only on points at the interface, we consider an extension of the velocity to the whole domain $Y$. This allows us to describe the cancer region boundary by a level-set approach. The interface is defined as the zero-level of the level-set function $\phi$ which satisfies the following transport equation:
\begin{equation}
\frac{\partial \phi}{\partial t}+V(m)\cdot\nabla\phi=0,\quad \text{ in } Y \times \tho{(0,T)}.
\label{eq:transportequation}
\end{equation}

\tho{
For later purposes, we introduce the notation 
\begin{equation}
\label{L0}
L_0(t) = \{x\in Y: \phi(x, t)=0\}
\end{equation}
for the zero level of the level set function that defines the interface $\partial \Omega(t)$.
}

A natural extension of the velocity is the constant continuation of the velocity at the boundary in normal direction \cite{Chopp:2009}. In Section \ref{sec:extension} more details about this point are given.

\subsection{Microscopic model component}
Due to the scale separation Assumption \ref{ass:scale separation} we can describe the dynamics of the \dtr{MDE} on a microscopic domain \dtr{$\varepsilon Y$} defined \dtr{and centred at each macroscopic interface point $x\in\partial \Omega(t)$ as follows.}
\label{sec:micro}
\dtr{As argued in \cite{trucu2013multiscale}, the tumour cells arriving within a certain radius $R_{m}$ from the interface location $x\in\partial \Omega(t)$ give rise to a source of MDE, which simply represents the collective secretion of matrix degrading enzymes by the cells from the outer proliferating rim of the tumour that get distributed during their macro-dynamics on $B:=\{\xi\in Y: \|\xi - x\|\leq R_m\}$. Thus, mathematically, this source can be formalised as: } 
\begin{equation}
\label{F}
F_{x\dtr{,t}}(c)(y):=\left\{ \begin{array}{l l}
\displaystyle\frac{1}{\vert B \vert}\int\limits_{B}c(\xi\dtr{,t})\, {\rm d}\xi &\quad y \in \eps Y \cap \Omega(t)\\
 0 & \quad \text{otherwise}.\end{array}\right.
\end{equation}

\dtr{Therefore, in the presence of source \eqref{F}, the cross-interface micro-dynamics of MDE molecular distribution $m(y,\tau)$, which takes place on the micro-scale domain $\varepsilon Y$ over a micro-scale time range $(0,\Delta T)$, is governed by the following reaction diffusion equation}
\begin{align}
\begin{aligned}
\frac{\partial m}{\partial t}(y,\tau) &= D_2\Delta m(y,\tau)+F_{x\dtr{,t}}(c) &&\quad \text{in}\  \eps Y \times (0,\Delta T)  \\
m(y,0) &= 0 && \quad  \text{in}\ \eps Y\\
\frac{\partial m }{\partial n}(y,\tau) &=0 && \quad  \text{in}\ \partial \eps Y \times (0,\Delta T),
\end{aligned}
\label{eq:microprob}
\end{align}
with $\Delta T$ \dtr{representing here the micro-scale time perspective and serving also later on as \emph{natural time splitting step} between micro and macro-stages within the computational approach of the multiscale model}.

\dtr{Finally, we would like to remark here the special multiscale importance of the microscale source term given via the non-local operator $F_{x\dtr{,t}}(c)$ that is induced by the macro-dynamics and realises a \emph{top-down} link from macro- to micro-dynamics, enabling this way the entire interface micro-dynamics. }\\ 

\paragraph{Time splitting}
At this stage, the three parts of the model, i.e.\ the macroscopic, the transport and the microscopic components, are fully coupled and the coupling is given at the continuous level.

To define the velocity \dtr{$V$} used in equation \eqref{eq:transportequation} we consider a splitting in time of the overall coupling between the three model components.
Therefore, for a given macroscopic time interval $\Delta T$\dtr{, the pattern of peritumoural ECM degradation caused by the advancing fronts of MDE molecules (which are transported across the tumour interface in the immediate proximity within the appropriate microscale region) gives rise to a boundary velocity  that can be mathematically described by}%
 \begin{equation}
 \dtr{V}(m)=\frac{\cv}{\Delta T\, |\eps Y|}\,\int_0^{\Delta T}\int_{\eps Y}m\, \nabla m\, \dtr{\dy}{\rm d} \tau \label{eq:velocity}
 \end{equation}
where $|\eps Y| = \int_{\eps Y} 1\; \dt$ and $\cv$ is a tuning scaling factor, see Table \ref{tab:num setting}. \dtr{Specifically, this form of $V$ is based on the following main considerations: }
\begin{itemize}
\item the term $\nabla m$ takes into consideration the assumption that the cancer boundary moves following the gradient with respect to the MDE;
\item \dtr{further}, by multiplying it by $m\dtr{(y,\tau)}$, we are taking into account the influence of the amount of enzymes \dtr{over their given gradient direction at each spatio-temporal micro-node $(y,\tau)$, enabling an appropriate weighting of its \emph{``strength"} (magnitude)};
\item \dtr{finally,} by considering the average contribution of \dtr{MDE microdynamics over $\varepsilon Y \times [0,\Delta T]$ by simply accounting upon the mean-value in time of the \emph{revolving weighted MDE gradient spatial direction} 
\[
 [0,\Delta T]\ni\tau \mapsto \frac{1}{|\eps Y|}\int_{\eps Y}m\, \nabla m\, \dtr{\dy},
 \]}\dtr{we ultimately obtain the definition of the velocity given in \eqref{eq:velocity}, where $V(m)$ is taken as being proportional to this spatio-temporal mean value, with proportionality constant $c_{vel}$}.
\end{itemize}

\begin{figure}[!h]
\begin{center}
        \scalebox{0.8}{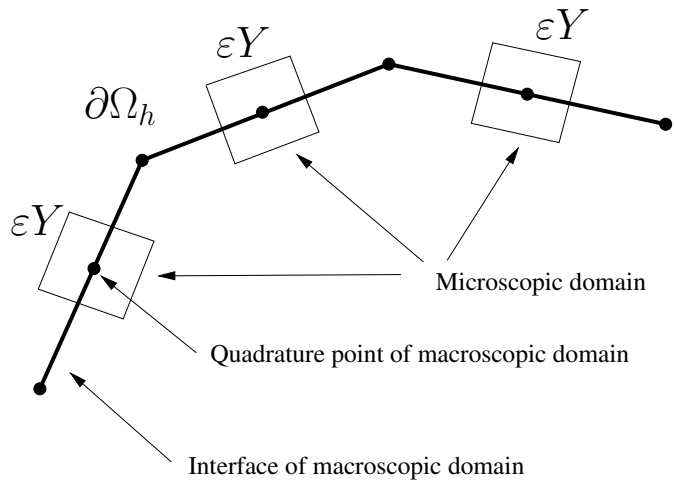}
\end{center}
\caption{Microdynamics sampling at the macroscopic boundary} 
\label{sampling}
\end{figure}

In the implementation we use a linearized interface. For each macroscopic cell that is cut by the interface, we consider quadrature points at the interface and assign a velocity vector to each of them, see Figure \ref{sampling} for a sketch of the interface.

\subsection{Definition of the computational microscopic problem}
To simplify the implementation of the numerical method we consider the microscopic problem in \dtr{a bundle of boundary micro}domain\dtr{s} $\eps Y$, where $Y$ is \dtr{the reference maximal tissue cube centred at origin}, \dtr{with} $\dtr{y:=}(y_1,y_2)$ \dtr{being} the standard local \dtr{microscale} reference system \dtr{within a given $\varepsilon Y$} \dtr{ and $\tau$ denoting always the time at micro-scale}.

As \dtr{will be} explained below in Section \ref{sec: linearized interface}, \dtr{in our finite element approach} the macroscopic \dtr{dynamics will be considered} on a \dtr{appropriately defined macroscopic} domain $\Omega_h(t)$ with a linearized boundary $\partial \Omega_h(t)$. 
The microscopic \dtr{dynamics} is then \dtr{explored} within a \dtr{microdomain $\varepsilon Y$ centred at a macro-scale boundary point $x\in\partial \Omega_h(t) $ and eventually appropriately rotated so that this is positioned} with two edges parallel to the linearized boundary (in direction $y_1$) and two edges orthogonal to it (in direction $y_2$) as shown in Figure \ref{sampling}. Th\dtr{is} orientation of the computational microscopic domain is well defined because it is centered in quadrature points at the boundary, which are never defined at the corners of the piece-wise linear boundary.
This simplifies the setting of the microscopic problem.
In fact, since we consider the linearized boundary $\partial \Omega_h(t)$, the right hand side $F_{x\dtr{,t}}$ in equation \eqref{eq:microprob} does not depend on $y_1$. In addition, since on the boundaries of the quadrilateral domain no flux conditions are prescribed, it follows that the solution is constant in $y_1$ direction. Therefore, we can consider the following simplified one-dimensional microscopic problem for the quantity $\overline m$, which is the integral of $m$ along $y_1$ (giving the amount of enzyme molecules per unit of length) 
\begin{align}
\begin{aligned}
\frac{\partial \overline{m}}{\partial \tau}(y_2,\tau) &= D_2\Delta \overline{m}(y_2,\tau)+\overline{F}_{x\dtr{,t}}(y_2) && \quad \text{in}\  (0,\eps)\times \tho{(0,\Delta T)}  \\
\overline{m}(y_2,0) &= 0 && \quad \text{in}\ (0,\eps)\\
\frac{\partial \overline{m}}{\partial n}(y_2,\dtr{\tau}) &= 0 && \quad \text{in}\ \partial (0,\eps)\times \tho{(0,\Delta T)}
\label{eq:microprob1D}
\end{aligned}
\end{align}
where $\overline{F}_{x\dtr{,t}}$ is $F_{x\dtr{,t}}$ integrated over $y_1$.
Since $F_{x\dtr{,t}}$ and $m$ do not depend on $y_1$, the solution $m$ of \eqref{eq:microprob} is the constant extension of the solution $\overline{m}$ of \eqref{eq:microprob1D} in $y_1$ direction.

We introduce now a scaling of the domain to the interval $(0,1)$ through the following tranformation
\begin{equation*}
y_2 =\eps z \quad \text{with }\ z\in (0,1), 
\end{equation*}
then we get after the rescaling the transformed system 
\begin{align}
\begin{aligned}
\frac{\partial \widehat{m}}{\partial t}(z,\dtr{\tau}) &= D_2 \eps^{-2} \Delta \widehat{m}(z,\dtr{\tau})+\widehat{F}_{x\tho{,t}}( z)&& \quad \text{ in }\  (0,1)\times \tho{(0,\Delta T)} \\
\widehat{m}(z,0)&=0 && \quad \text{ in }\ (0,1)\\
\frac{\partial \widehat{m}}{\partial n}(z,\dtr{\tau})&=0 && \quad \text{ in }\ \partial (0,1)\times \tho{(0,\Delta T)}
\label{eq:microprob1Dad}
\end{aligned}
\end{align}
with
\begin{equation}
\label{widehatF}
\widehat{F}_{x\dtr{,t}}(z):=\left\{ \begin{array}{l l}
\displaystyle\frac{1}{\vert B\vert}\int\limits_{B}c(\xi\dtr{,t})\, {\rm d}\xi  & z \in [0,1/2]\\
 0 & \text{otherwise},\end{array}\right.
\end{equation}
note that the coordinate $\xi$ is a macroscopic quantity.
Notice furthermore that a solution of \eqref{eq:microprob1D} is a solution of \eqref{eq:microprob1Dad} by $\widehat{m}(z\dtr{,\tau})=\overline{m}(\eps z\dtr{,\tau})$.
\begin{remark}[Limit $\eps \rightarrow 0$]
In case of $\eps \rightarrow 0$ we have in \eqref{eq:microprob1Dad} a large diffusion coefficient, therefore a fast redistribution process of the solution occurs, leading to negligible spatial variations of the solution. The only relevant parameter of the problem at the limit becomes the time. The limit problem becomes an ordinary differential equation (ODE).
Even if we consider scale separation in this model, we do not consider the limiting case $\eps \rightarrow 0$. In that case the velocity has to be defined in a different way since the term $\nabla m$ becomes the zero vector. The parameter $\eps$ in our model has always a finite value bounded below $\eps \geq \eps\dtr{_{_{cell}}} > 0$, \dtr{where $\eps_{_{cell}}$ is assumed here to be a minimal microscale size of the order of a cell-length.}
\end{remark}

\dtr{Thus, using \eqref{eq:microprob1Dad}, we obtain that the velocity defined by} problem \eqref{eq:microprob} via equation \eqref{eq:velocity} \dtr{can be further expressed as}
 \begin{align}
\nonumber%label{velocity}
V(m) = \frac{\tho{\cv}}{\Delta T \eps^2} \int_{[0,\eps]^2}\int_0^{\Delta T}m\nabla_y m\, {\rm d}\dtr{\tau}\, {\rm d}y&=\frac{\tho{\cv}}{\Delta T\eps^2}\eps \int_0^\eps\int_0^{\Delta T}m\, \nabla_y m\, \dy\, {\rm d}\dtr{\tau}\\
\label{scaled velocity}
&=\frac{\tho{\cv}}{\Delta T\, \eps}\int_0^1\int_0^{\Delta T}\widehat{m}\, \nabla_z\widehat{m}\, \dz\, {\rm d}\dtr{\tau}
 \end{align}
where $\widehat{m}$ indicates the transformed function on the \dtr{reference} domain \dtr{$Y$}.

\section{Weak formulation and discretized model}
\label{weak form}
To describe the model in the setting needed for the FEM we introduce the following weak formulation.
\subsection{Weak formulation}
We use the notation $(\cdot, \cdot)$ to define the usual $L^2$ scalar product of Lebesgue square integrable functions. The space $H^1$ is the Hilbert space of square integrable functions with square integrable (weak) first derivative and $H^*$ is its dual space, i.e.\ the space of bounded linear functional on $H^1$. Furthermore, we use Bochner spaces like $U = \{u \in L^2(0,T;H^1): \partial_t u \in L^2(0,T; H^*)\}$ to introduce the weak formulation of each subproblem.
In particular, we consider the functional space \[U_T=\{\phi \in L^2(0,T;H^1(Y)): \partial_t \phi \in L^2(0,T; H^*(Y))\}\] for the transport component, the space \[U_M=\{m \in L^2(0,T;H^1(\Omega(t))): \partial_t m \in L^2(0,T; H^*(\Omega(t)))\}\] for the macroscopic component and \[U_m = \{c \in L^2(0,T;H^1((0,1))): \partial_t c \in L^2(0,T; H^*((0,1)))\}\] for the scaled microscopic component. The weak formulation problem of the macroscopic part is: find the pair $(c,v)\in U_M \times U_M$ such that it satisfies for almost all $t \in (0,T)$ 
\begin{alignat*}{2}
\big( \dvar{c}, \varphi \big) + \big(D_1 \nabla c, \nabla \varphi \big) - \big( \eta c \nabla v, \nabla \varphi \big) + \big( \mu_1\tho{(v)}\, c (1-c-v), \varphi \big) &= 0 &&\, \forall \varphi \in H^1(\Omega(t))\\
\big( \dvar{v}, \varphi  \big) + \big( \alpha c\,v, \varphi \big) + \big( \mu_2 (1-c-v), \varphi \big) &= 0 &&\, \forall \varphi \in H^1(\Omega(t))\\
c(x, 0) &= c_0 && \, \text{ in } \Omega(0)\\
v(x, 0) &= v_0 && \, \text{ in } \Omega(0)
\end{alignat*}
Note that the cancer cells $c_0$ and ECM distributions $v_0$ are defined on the larger \dtr{maximal} domain $Y$ and that this formulation implies the {\em natural} zero flux condition for the cancer cells and ECM, i.e.\ $\partial_n c = \partial_n v =  0 \text{ on } \partial \Omega(t)$.

The weak formulation of the transport equation is: find $\phi \in U_T$ such that for a.a.\ $t \in (0,T)$ it satisfies
\begin{alignat*}{2}
\big( \dvar{\phi}, \varphi \big) + \big( V(m) \cdot \nabla \phi, \varphi \big) &= 0 && \quad \forall \varphi \in H^1(Y)\\
\phi(x, 0) &= \phi_0 && \quad \text{ in } Y,
\end{alignat*}
with $\phi_0$ being the level set function at the initial time. 
Using the notation \[\varphi^+:=\displaystyle \max\{\varphi, 0\}\] to define the positive part of a function $\varphi$ in its domain of definition, we have that $\text{supp}(\phi_0^+)$ describes the initial support region of the cancer cells.

The weak formulation of the microscopic problem is: find $\mh \in U_m$ such that for a.a.\ $\tau \in (0,\Delta T)$ it satisfies
\begin{alignat*}{2}
\big( \dvarp{\mh}, \varphi \big) + \big(D_2 \eps^{-2} \nabla \mh, \nabla \varphi \big) &= \big( \widehat F_{x\dtr{,t}}, \varphi \big) && \quad \forall \varphi \in H^1((0,1))\\
\mh(z, 0) &= 0 && \quad  \text{ in } (0,1),
\end{alignat*}
where $\widehat F_{x\dtr{,t}}$ is defined as in \eqref{widehatF} and the natural condition $\partial_n \mh = 0$ is implicitly defined.
\subsection{Discretization}
The model is first discretized in time by the implicit Euler method and then discretized in space by the FEM. 
The discretized system is defined on a regular mesh $\T$ composed of quadrilateral cells $K \in \T$ of the same dimension. This discretization has in a natural way the property of shape regularity, \tho{i.e.\ there exists a constant $\kappa$ such that $\displaystyle \max\limits_{K\in \T} \frac{h_K}{\rho_T} \leq \kappa$ for all elements of the mesh}, where $h_K$ and $\rho_K$ are respectively the cell diameter and the diameter of the largest ball inscribed into $K$. Furthermore it has the advantage that the mesh can be generated starting from an initial quadrilateral that is refined successively to reach a given cell diameter.

Since the macroscopic domain $\Omega(t)$ is time-dependent,
the discrete space domain would need to be remeshed at every time step, if a fitted FEM formulation is used.
In case of large deformations, the procedure of remeshing
has to deal with the possible loss of shape regularity of the mesh. 

To avoid these complications related to remeshing, we use an unfitted approach by using so called {\em cut-cells}. These are a special realization of the FEM as described below.
In particular, these are finite elements with shape functions with a support on a subdomain of the cells that is defined by the intersection of the interface with the cells. 

\begin{figure}[!h]
\begin{center}
        \scalebox{0.4}{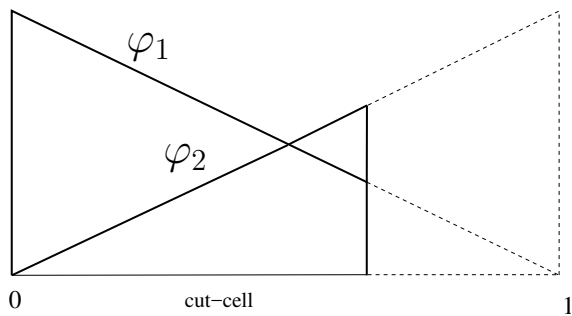}
\end{center}
\caption{Shape functions on a unit cut-cell}
\label{fig:cut-cells}
\end{figure}

Let us consider the space of bi-linear polynomials $Q_1$ defined on a unit cell $\hat K=[0,1]^2$, i.e.\
\[Q_1=\text{span}(1,x,y,xy)\]
\tho{and the space of linear functions
\[P_1=\text{span}(1,x)\]
defined on the unit one-dimensional cell $\hat K=[0,1]$.}
The finite element space is defined as
\begin{align*}
\VMh(t) = \{ u \in C(\Omega(t)): &u_{|K}\circ T_K \in Q_1 \text{ if } K\cap \partial \Omega(t) = \emptyset;\\ &u_{|K\cap \Omega(t)} = \psi_{|\Omega(t)},\, \psi \circ T_K \in Q_1 \text{ if } K\cap \partial \Omega(t) \neq \emptyset  \},
\end{align*}
where $T_K$ is a bijective transformation from the unit cell $\hat K$ to the physical cell $K$. The functions $v$ are called pull-back functions and are the transformations to the real coordinates of the polynomial functions defined on the unit cells. In our case, we consider only a translation and a scaling of the unit cells. Since the mesh is non-fitted, the functions $v$ (in case of cut-cells) are defined only on the portion of cell that is intersected by $\Omega(t)$.
On the portion of the cell that lies outside of $\Omega(t)$ the shape functions need not to be defined. \tho{In Figure \ref{fig:cut-cells} a one dimensional sketch is shown, where one can see how the restriction of the shape functions $\varphi_1$ and $\varphi_2$ on the cut-cell is defined. Note that the shape functions are not modified, we just use a restriction of them on the part of the cell that belongs to the domain. The Lagrangian formulation remains the conventional one for continuous finite elements. The degrees of freedom are defined on the cell nodes as usual. 
}
Therefore, although the cells are arbitrarily cut by the boundary $\partial \Omega(t)$, the convergence rate of this finite element formulation is not affected by the position of the cut, leading to optimal discretizations. This result can be directly derived by the convergence proof shown in \cite{Hansbo:2002}. 

For the transport and microscopic problems we use the following spaces
\begin{align*}
\VTh = \{ u \in C(Y): u_{|K} \circ T_K \in Q_1 \}
\end{align*}
and
\begin{align*}
\Vmh = \{ u \in C((0,1)): u_{|K} \circ T_K \in P_1 \}.
\end{align*}
\paragraph{Approximation of the moving domain}
For the time discretization of the macroscopic problem the variation of the domain in time should be taken into account. 
One possible formulation of the problem would be to consider a reference domain $\Omega(t_{0})$ for a given $t_0$ in each time step and to use a (time dependent) mapping to transform the solution from the domain $\Omega(t_n)$ to the reference domain $\Omega(t_{0})$ \tho{and then to the domain $\Omega(t_{n+1})$}. 
However, if the time step is small enough \tho{the combination of these two transformations} can be approximated with the identity and instead an extension of the solution from the old domain $\Omega(t_n)$ to the new domain $\Omega(t_{n+1})$ could be used. This procedure introduces an approximation error that for small enough time steps can be neglected in comparison to other sources of error (\dtr{such as} space-time discretization of the solution, \dtr{or} splitting error etc.). 
\tho{Due to the complexity of the formulation of the presented approach we leave the investigation of other moving domain formulations for a forthcoming work}.
\begin{figure}[!h]
\begin{center}
        \scalebox{0.35}{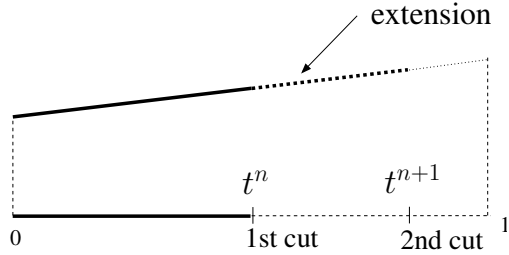}
\end{center}
\caption{One dimensional sketch of the extension of the macroscopic solution in case of a cell cut twice}
\label{fig:extension}
\end{figure}
\begin{figure}[!h]
\begin{center}
        \scalebox{0.35}{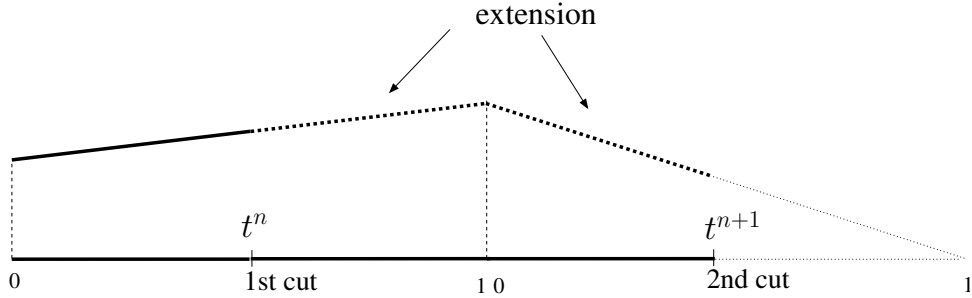}
\end{center}
\caption{One dimensional sketch of the extension of the macroscopic solution in case the interface cuts two neighbour cells at $t^n$ and $t^{n+1}$.}
\label{fig:extension2}
\end{figure}

Following the above construction, the fully discrete formulation of the macroscopic problem becomes
\begin{subequations}
\begin{align}
\begin{split}
&\big( c_h^{n+1}, \varphi \big)_{\Omega(t_{n+1})} + k \big[ \big(D_1 \nabla c_h^{n+1}, \nabla \varphi \big)_{\Omega(t_{n+1})} - \big( \eta\, c_h^{n+1} \nabla v_h, \nabla \varphi \big)_{\Omega(t_{n+1})}\\[2mm]
&+ \big( \mu_1\tho{(v_h^{n+1})}\, c_h^{n+1} (1-c_h^{n+1}-v_h^{n+1}), \varphi \big)_{\Omega(t_{n+1})} \big] = \big( \tilde c_h^{n}, \varphi \big)_{\Omega(t_{n+1})}, \qquad \forall \varphi \in \VMh(t^{n+1}),
\end{split}\label{discrete u} \\[5mm]
\begin{split}
&\big( v_h^{n+1}, \varphi  \big)_{\Omega(t_{n+1})} + k \big[ \big( \alpha c_h^{n+1}\,v_h^{n+1}, \varphi \big)_{\Omega(t_{n+1})} + \big( \mu_2\, (1-c^{n+1}_h-v_h^{n+1}), \varphi \big)_{\Omega(t_{n+1})} \big] \\[2mm]
&= \big( \tilde v_h^{n}, \varphi  \big)_{\Omega(t_{n+1})}, \qquad \forall \varphi \in \VMh(t^{n+1}),
\end{split}\label{discrete v}
\end{align}\label{discrete:macro-sys}
\end{subequations}
\noindent where $\tilde c_h^n$ and $\tilde v_h^n$ are the extensions from $\Omega(t_n)$ to $\Omega(t_{n+1})$. In fact, the two components $c_h(t_n)$ and $v_h(t_n)$ are defined only in $\Omega(t_n)$. Therefore an extension in the region $\Omega(t_{n+1}) \setminus \Omega({t_n})$ needs to be defined. We have chosen a continuous extension using the prescribed values $c_0(x)$ and $v_0(x)$ for $x \in Y$. In particular, we have considered two cases: case (i) the cell where we need to define the extension is cut at time $t^n$ and at time $t^{n+1}$, see Figure \ref{fig:extension} and case (ii) the cell is cut at time $t^n$ and uncut at time $t^{n+1}$, see Figure \ref{fig:extension2}. In case (ii) the cut goes to the neighbour cell at time $t^{n+1}$. In case (i) both components are extended up to the new cut using the values of all degrees of freedom of the considered cell (also those \tho{lying outside} the domain $\Omega(t)$) with a bilinear nodal interpolation. Note that the bilinear nodal interpolation is justified only within the domain $\Omega(t)$ in the cut-cell formulation of the problem, because the integrals in the weak formulation are computed only in the inner part of the cells. In this sense, we are ``extrapolating'' the values $c_h$ and $v_h$ outside the region of validity of the finite element interpolation. We have observed that this ``extrapolation'' gives good results only for cut-cells that are not too small. In case of small cut-cells, the extrapolation leads to wrong values that introduce instabilities, visible as large peaks in the solution, that destroy the convergence of the method. We have set a threshold of 1\% on the volume to be considered for extrapolation. Cells, whose volume is cut by 99\%, are eliminated from the active mesh. We are currently studying the reason for this behavior of small cut-cells and let the rigorous discussion of this issue for a subsequent work since it goes beyond the scope of this paper. 

The system of equations is solved with the following initial conditions
\begin{align*}
c_h^0 &= c_0  \quad \text{ in } \Omega(0),\\
v_h^0 &= v_0  \quad \text{ in } \Omega(0).
\end{align*}
The transport equation is defined on $Y$. Since this is a hyperbolic equation, a suitable discretization is needed. We choose the stream-line diffusion approach for its easy implementation and good performance. We have used an artificial diffusion in the stream-line direction scaled with a parameter $\delta > 0$ whose value can be found in Table \ref{tab:num setting}. On general meshes, the stream-line diffusion stabilized formulation of the transport equation converges in the $L^2$ norm with the rate $h^{3/2}$ with bilinear finite elements. Nevertheless, several authors show an optimal convergence rate of $h^2$ on regular meshes, see for example \cite{Zhou:1997}.

The dynamics of the cancer cells and ECM are solely defined by the velocity at the boundary and the initial distributions $c(x,0)$ and $v(x,0)$.
\tho{The discrete transport equation is
\begin{equation}
\begin{array}{rll}\label{discrete phi}
\big( \phi_h^{n+1}, \varphi \big)_{Y} + k \big( V_h^n\cdot \nabla \phi_h^{n+1}, \varphi + \delta\, (V_h^n\cdot \nabla \varphi) \big)_{Y} &= \big( \phi_h^n, \varphi \big)_{Y} & \, \forall \varphi \in \VTh\\[2mm]
\phi_h^0 &= \phi_0 & \, \forall x \in\, Y,
\end{array}
\end{equation}
where $\phi_0$ is the initial level set function, $k$ the time step and $V_h^n$ is the discrete velocity defined as
\begin{equation}
\label{discrete velocity}
V_h^n:= \frac{\cv}{\Delta T\, \eps} I_x \big(I_\tau\big(\mhhn\, \nabla \mhhn\big)\big),
\end{equation}
where $I_x$ and $I_\tau$ are two quadrature formulas for the approximation of the integral in space and time (see expression \eqref{scaled velocity}) and $\mhhn$ is the discrete solution of the microscopic problem as defined below.
}

Since we assume scale separation, in each point of the macroscopic boundary $\partial \Omega(t)$ we need to solve a microscopic problem that defines the local velocity.
In the discrete version, we define the microscopic problem in a finite number of points at the interface and we discuss later the issue of how to use these point wise defined velocities to solve the transport problem.
\tho{
The weak formulation of the microscopic problem reads:
\begin{equation}
\begin{array}{rll}\label{discrete m}
\big( \mhhn^{l+1}, \varphi \big) + k \big(D_2 \eps^{-2} \nabla \mhhn^{l+1}, \nabla \varphi \big) &= \big( \mhhn^l, \varphi \big) + \big( \widehat F_{x,n,h}, \varphi \big) & \quad \forall \varphi \in \Vmh\\[2mm]
\mhhn^0 &= 0 & \quad \text{ in } (0,1),
\end{array}
\end{equation}
where we have used the notation $\mhhn^l$ to indicate the discrete microscopic solution for the macroscopic step $n$ (note that the right hand side depends on the macroscopic solution at time $t_n$), with $l$ being the time step of the time variable $\tau$, i.e.\ $\mhhn^l = \mhhn(t_l)$.
}
The term $\widehat F_{x,\dtr{n,}h}$ is an approximation of $\widehat F_{\dtr{x,t_{n}}}$ in which $c$ is substituted by its discrete counterpart and the integral over $B$ is approximated by a quadrature rule
\begin{equation}
\label{widehatFh}
\widehat{F}_{x,\dtr{n,}h}(z):=\left\{ \begin{array}{l l}
\displaystyle\frac{I_{B}(c^{\dtr{n}}_h)}{I_{B}(1)}  & z \in [0,1/2]\\
 0 & \text{otherwise},\end{array}\right.
\end{equation}
where $I_{B}(\cdot)$ is a quadrature rule that approximates the integral of the argument over $B$
\begin{align*}
I_{B}(f) \approx \int_{B} f(\xi)\, {\rm d} \xi.
\end{align*}
\subsection{Approximation of the interface and cut-cells}
\label{sec: linearized interface}
\begin{figure}[!h]
\begin{center}
\scalebox{0.2}{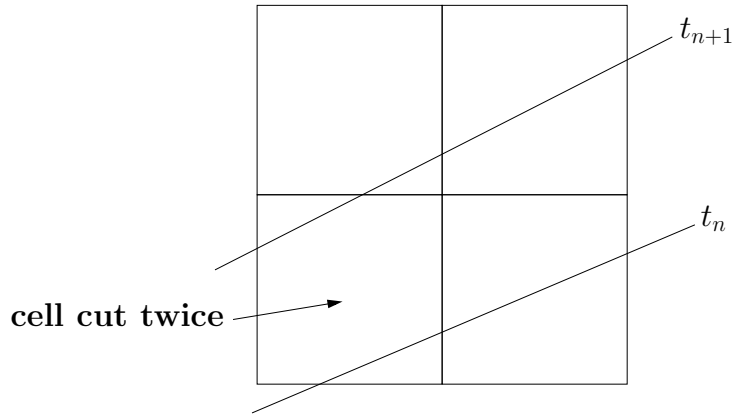}
\end{center}
\caption{Cell cut twice by the interface at two subsequent time steps.}
\label{twocuts}
\end{figure}

As introduced previously we discretize the problems in time by a time step method. Therefore, in the semidiscrete formulation we have terms that are defined at time $t=t_{n+1}$ and terms defined at time $t=t_n$. Since the domain is time dependent, the integrals of these terms are defined on different domains. Therefore, as explained above, in each time step we need to consider two configurations defined by the position of the boundary $\partial \Omega(t)$ in two subsequent time steps. In particular, we have to consider the case in which the boundary cuts the same finite element cell in both time steps, see in Figure \ref{twocuts} the cell at the bottom left and the case in which the boundary cuts one cell at time $t_n$ and it goes over to the neighbour cells at time $t_{n+1}$ leaving the previous cell uncut at time $t_{n+1}$, see in Figure \ref{twocuts} the cell at bottom right.

\begin{figure}[ht]
\centering
\subfloat[\label{subfig.lls}]{
\includegraphics[width=.4\textwidth]{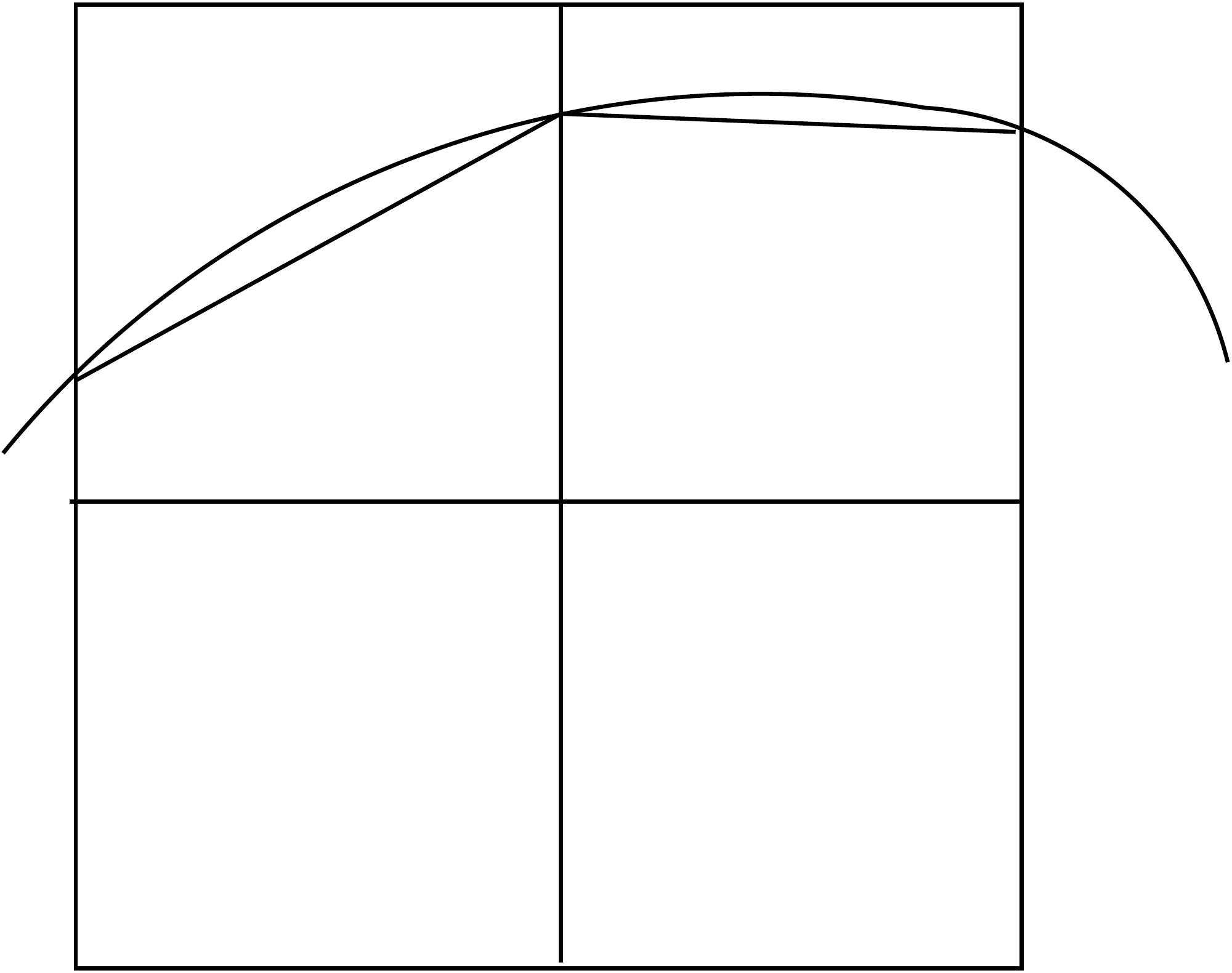}}
\subfloat[\label{subfig.mesh}]{
\includegraphics[width=.4\textwidth]{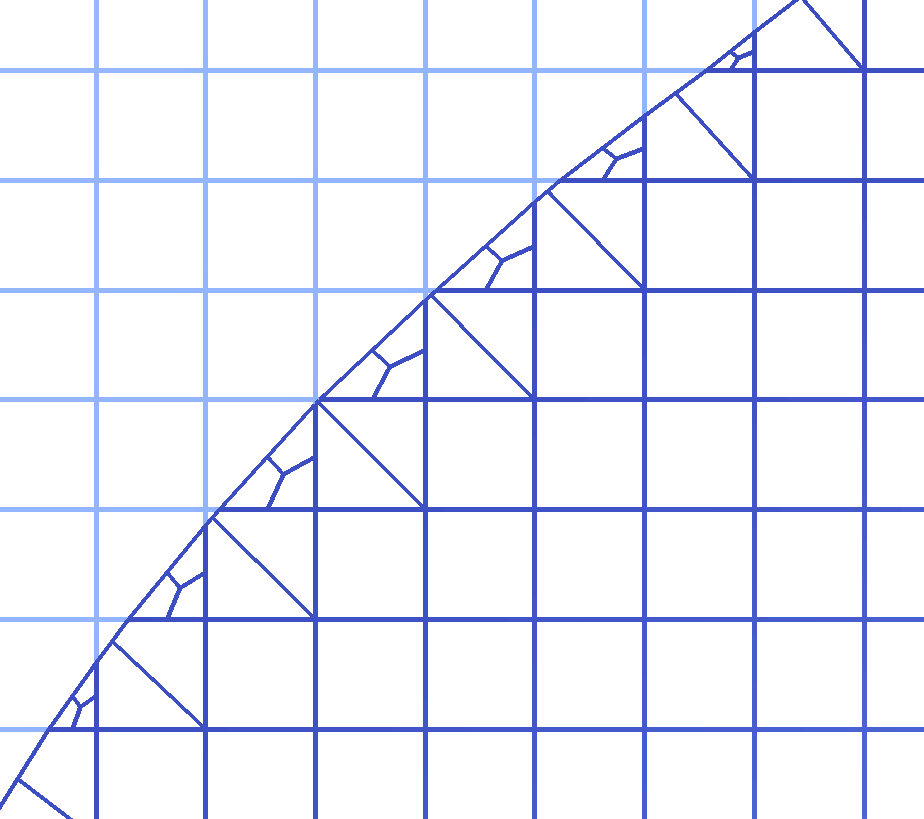}}
\label{linearized interface}
\caption{Left: Sketch of the linearized zero level $L_{0,h}$. Right: Portion of the mesh with the actual zero level and the cut cells.}
\end{figure}

For the discretized version of the system of equations, we consider the linearized domain $\Omega_h(t)$, which is defined by the piece-wise linear boundary
\begin{align}
\label{suppphi+}
\partial\Omega_h(t) := L_{0,h}(t),
\end{align}
where the linearized zero level $L_{0,h}(t)$ is defined by the polygonal line that connects all intersections of the zero level $L_0(t)$ with the mesh cells boundaries $\partial K$ as shown in Figure \ref{subfig.lls}. In Figure \ref{subfig.mesh} a part of the actual mesh is shown.
Furthermore, a part of the actual level zero isoline and the cut-cells are shown. Note that the irregular cells are not finite element cells, but are shown here only for visualization purpose.
By the linear approximation of the boundary we introduce an integration error in the solution of the weak formulation. It can be shown that this error converges to zero with the order of the discretization error. Therefore, the integration error can be interpreted as a perturbation of the underlying Galerkin method that does not change the convergence rate of the method.

Using the linearized boundary $\partial \Omega_h(t)$ we can apply the quadrature rule described in \cite{CarraroWetterauer:2016} to integrate the terms of the model on cut-cells with a single cut. Furthermore, if a cell is cut twice, i.e.\ at time $t_n$ and at time $t_{n+1}$, we apply the previous quadrature rule recursively.
\subsection{Approximation of the nonlocal term}

The nonlocal term \eqref{widehatFh} is approximated by a quadrature rule. We use a second level-set function to define the distance from the macroscopic point $x$ that is used to define the domain of integration $B$. Also in this case we introduce a piece-wise linear approximation of this second level-set function and we apply recursively the quadrature rule on the neighboring cells.

\subsection{Extension of the velocity}
\label{sec:extension}
The velocity is determined using formula \eqref{scaled velocity}.
It is computed at the macroscopic point $x$ on the linearized boundary $\partial \Omega_h(t)$ and then extended to the rest of the domain. 
In this work, we consider only one point $x$ per cell. This is taken at the midpoint of the segment of the interface $L_{0,h}$ that intersects the cell, as shown in Figure \ref{sampling}.

We set the velocity computed at this point $x$ to all cells which center lies at the closest distance from $x$. Therefore, the velocity is approximated as a piecewise constant function. For cells that lie in the cancer region, i.e.\ $K\cap \Omega_h(t)\neq 0$, at a distance larger than a prescribed radius of influence $\rho$ (see Table \ref{tab:num setting}) we set velocity zero. This procedure avoids the transport of numerical pollution from the center of the domain due to the singularity of the level set in the point that we take as reference to compute the distance function.

This definition of the velocity extension can lead to regularity problems in the transport of the interface if two parts of the boundary approach each other. This happens because the velocity of cells that lie at the same distance from the two approaching boundary parts is not well defined.
This is a typical problem in level set approaches that can be overcome using a fast marching method \cite{Adalsteinsson:1999}. Fast marching algorithms are often used in the context of level set approaches. Important issues in this context are the initialization of the level set function \cite{Chopp:2001}, its re-initialization as a signed distance function \cite{McCaslin:2014_2, McCaslin:2014_1, Gremaud:2007} and the possibility to update the extension of the velocity using an accurate characteristics reconstruction in case the velocity changes rapidly in time \cite{Chopp:2009}. 

\tho{
\subsection{Solution process}
We sketch the overall solution process underlying the coupling between the different parts of the model.
\begin{algorithm}[h]
\begin{algorithmic}[1]
\STATE Set $n=0$ and choose the splitting time step $\Delta T$
\STATE Set $\phi_0(x)$, $c_0(x)$ and $v_0(x)$ in $Y$
\STATE Define $L_0(t^n)$ as in \eqref{L0} and linearize it to get $\Omega_h(t^n)$
\STATE Solve macroscopic part \eqref{discrete:macro-sys} for $(x, t) \in \Omega_h(t^n)\times (t^n,t^n+\Delta T)$
\STATE Compute $\widehat{F}_{x,n,h}(z)$, see \eqref{widehatFh}
\STATE Solve microscopic part \eqref{discrete m} for $(x,\tau) \in \varepsilon Y \times (0, \Delta T)$
\STATE Compute velocity $V_h^n$, see \eqref{discrete velocity}
\STATE Extend velocity on all $Y$
\STATE Solve transport problem \eqref{eq:transportequation} for $(x, t) \in Y\times (t^n,t^n+\Delta T)$
\IF {$t^n+\Delta T = T$} 
\STATE \textbf{stop}
\ELSE  
\STATE Set $t^n= t^n + \Delta T$
\ENDIF
\STATE \textbf{goto} 3
\end{algorithmic}
\caption{Overall solution process}
\label{algo:solution}
\end{algorithm}
}

The macroscopic system is solved with an implicit \dtr{E}uler scheme. At each time step a nonlinear system of the type (\ref{discrete u}-\ref{discrete v}) has to be solved. We use an exact Jacobian and no damping for the Newton method, which converges generally in 2 steps to an accuracy lower than $10^{-6}$. The linear system arising in each Ne\dtr{w}ton step is solved by a direct solver. 

The system \eqref{discrete phi} is a linear system and is computationally much cheaper than the macroscopic problem. A direct solver is used in every time step.

Finally, the microscopic problem \eqref{discrete m} is a linear one dimensional parabolic problem solved with an implicit \dtr{E}uler method and a direct solver in each time step.

\begin{figure}[!h]
\centering
\subfloat[\label{subfig.}]{
	\includegraphics[width=.5\textwidth]{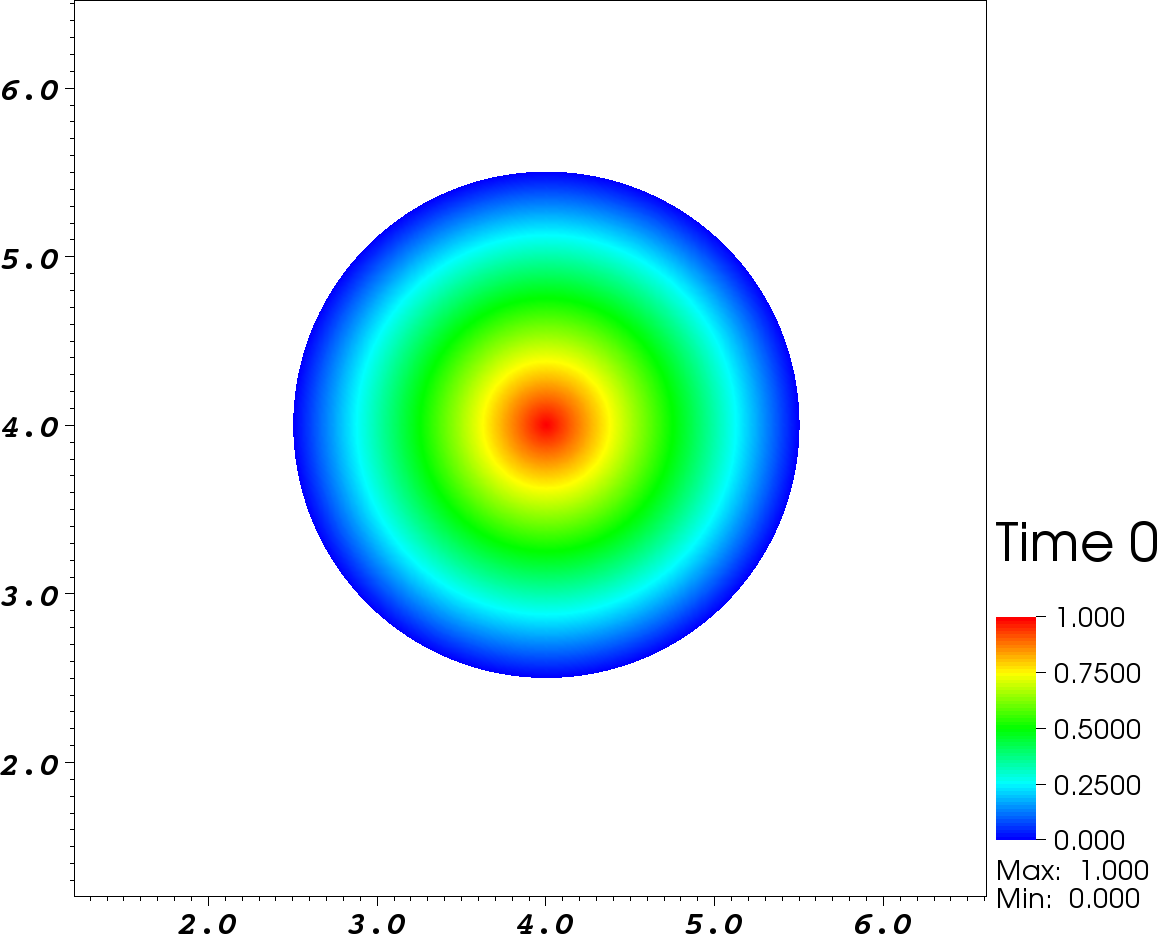}}
\subfloat[\label{subfig.}]{
	\includegraphics[width=.5\textwidth]{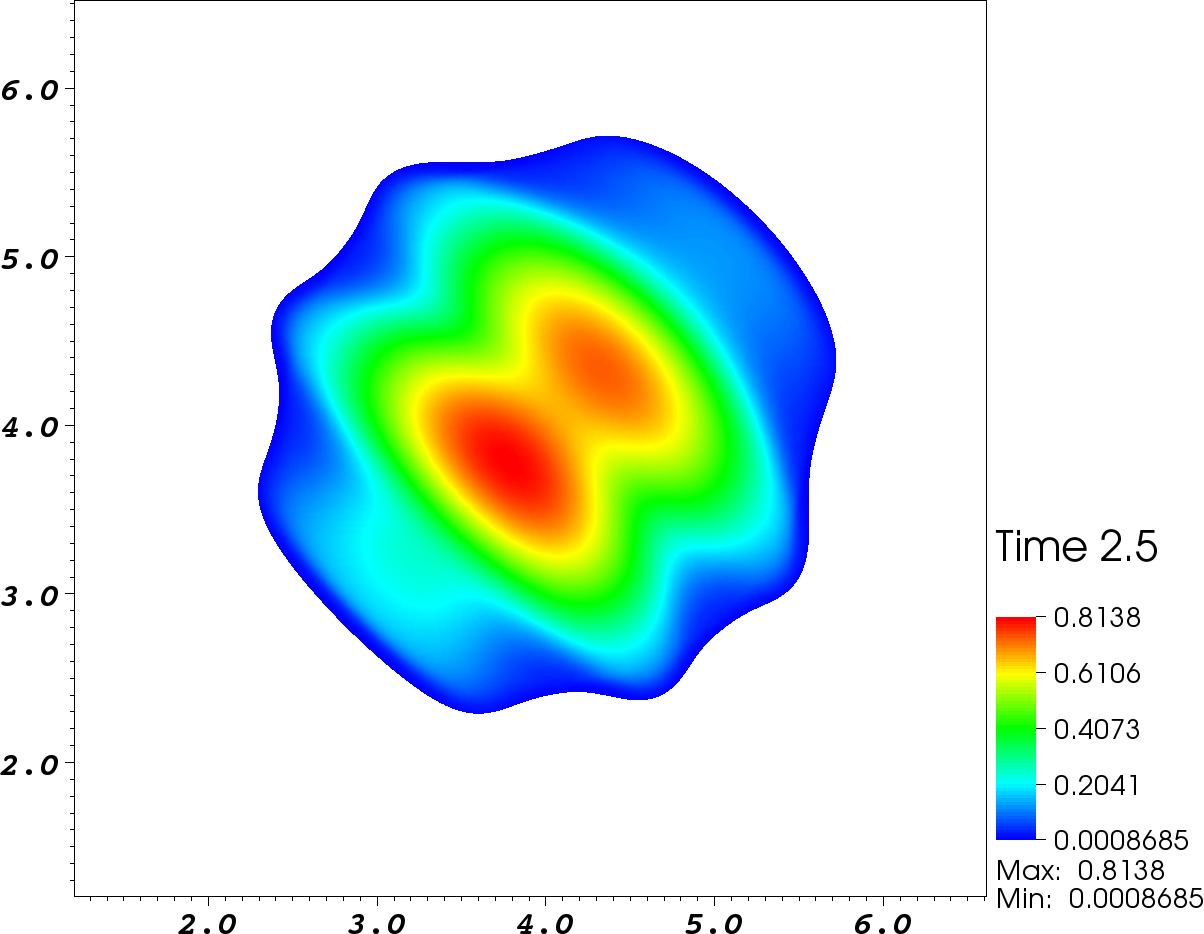}}\\
\subfloat[\label{subfig.}]{
	\includegraphics[width=.5\textwidth]{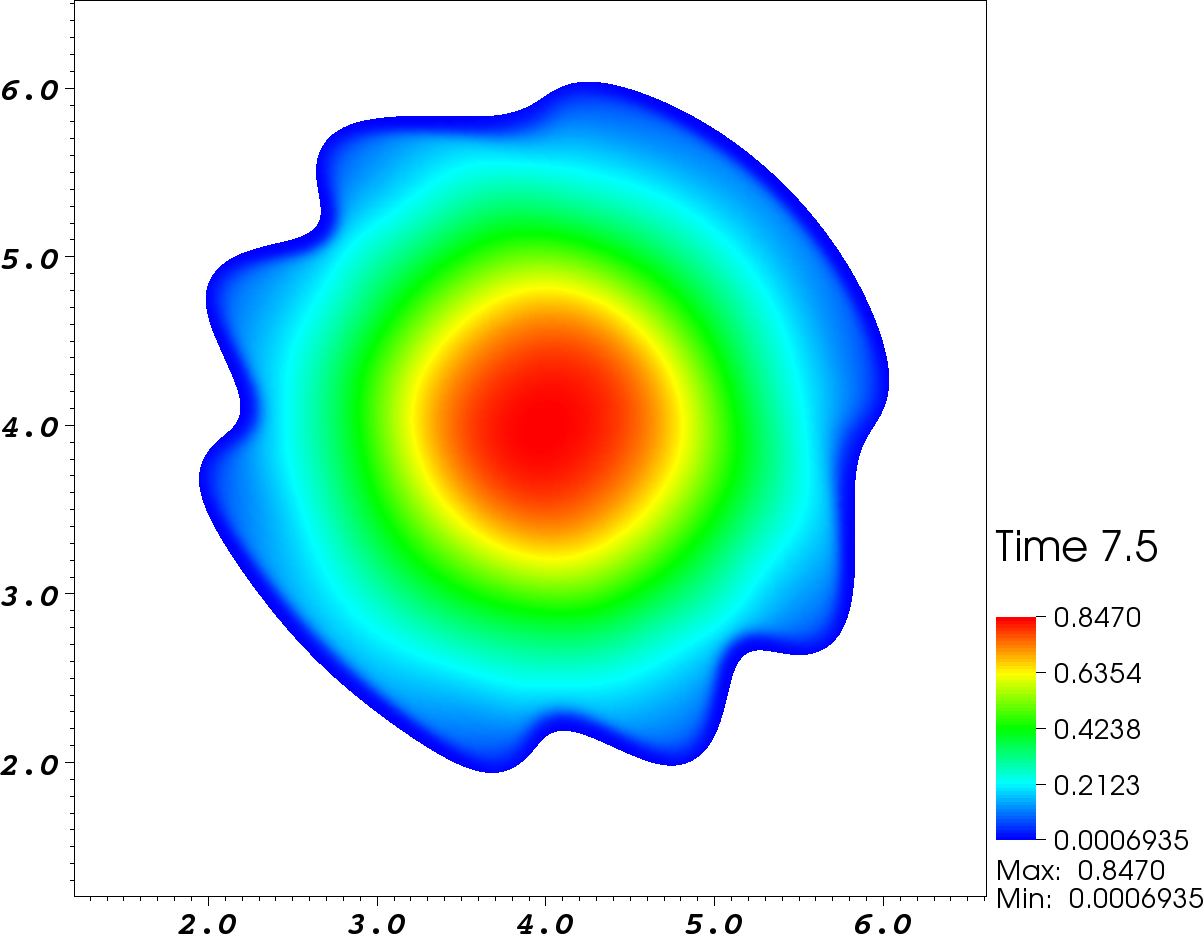}}
\subfloat[\label{subfig.}]{
	\includegraphics[width=.5\textwidth]{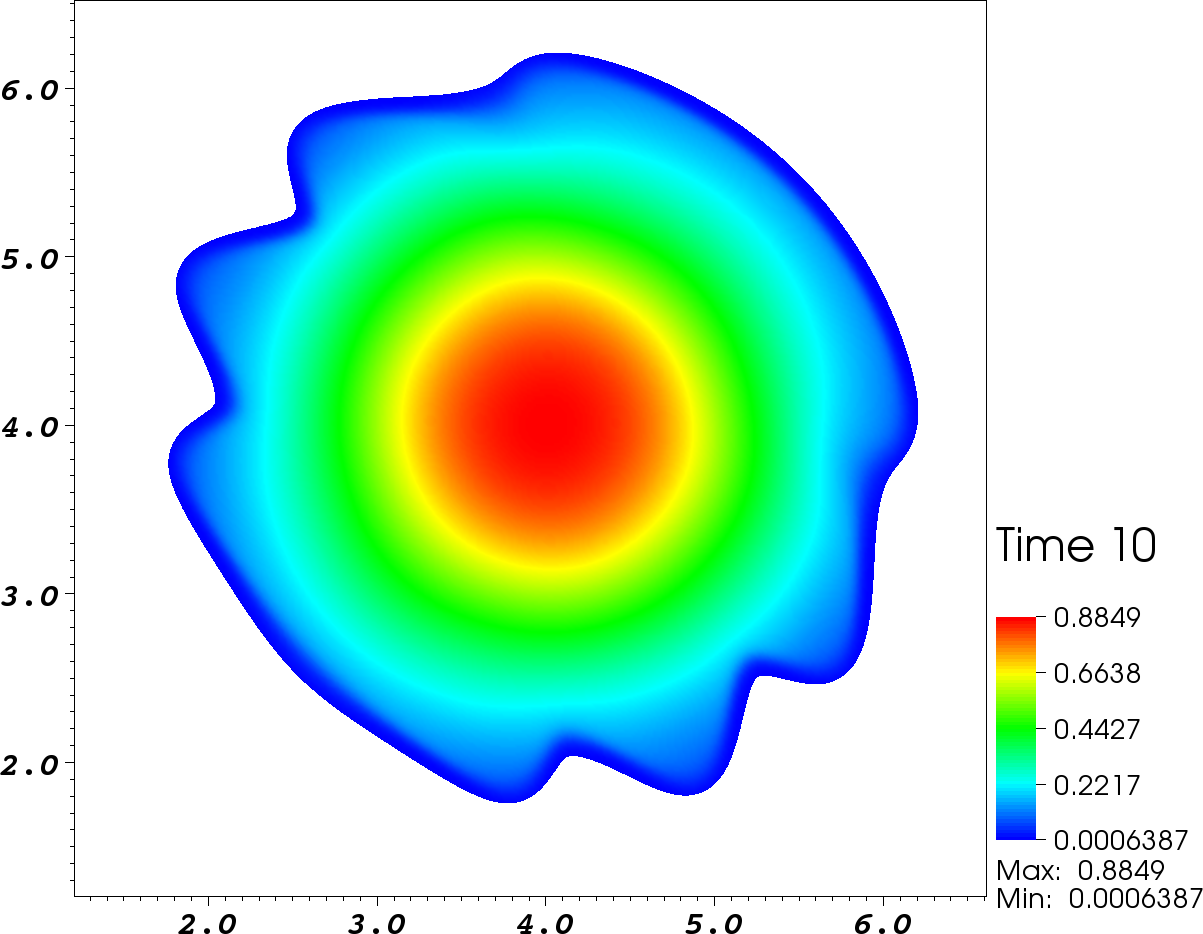}}
\caption{Cancer distributions at time\dtr{s: (a) $t=0$; (b) $t=2.5$; (c) $t=7.5$; and (d) $t=10$.}}
\label{distributions cancer}
\end{figure}

\begin{figure}[!h]
\begin{center}
\subfloat[\label{subfig.}]{
	\includegraphics[width=.5\textwidth]{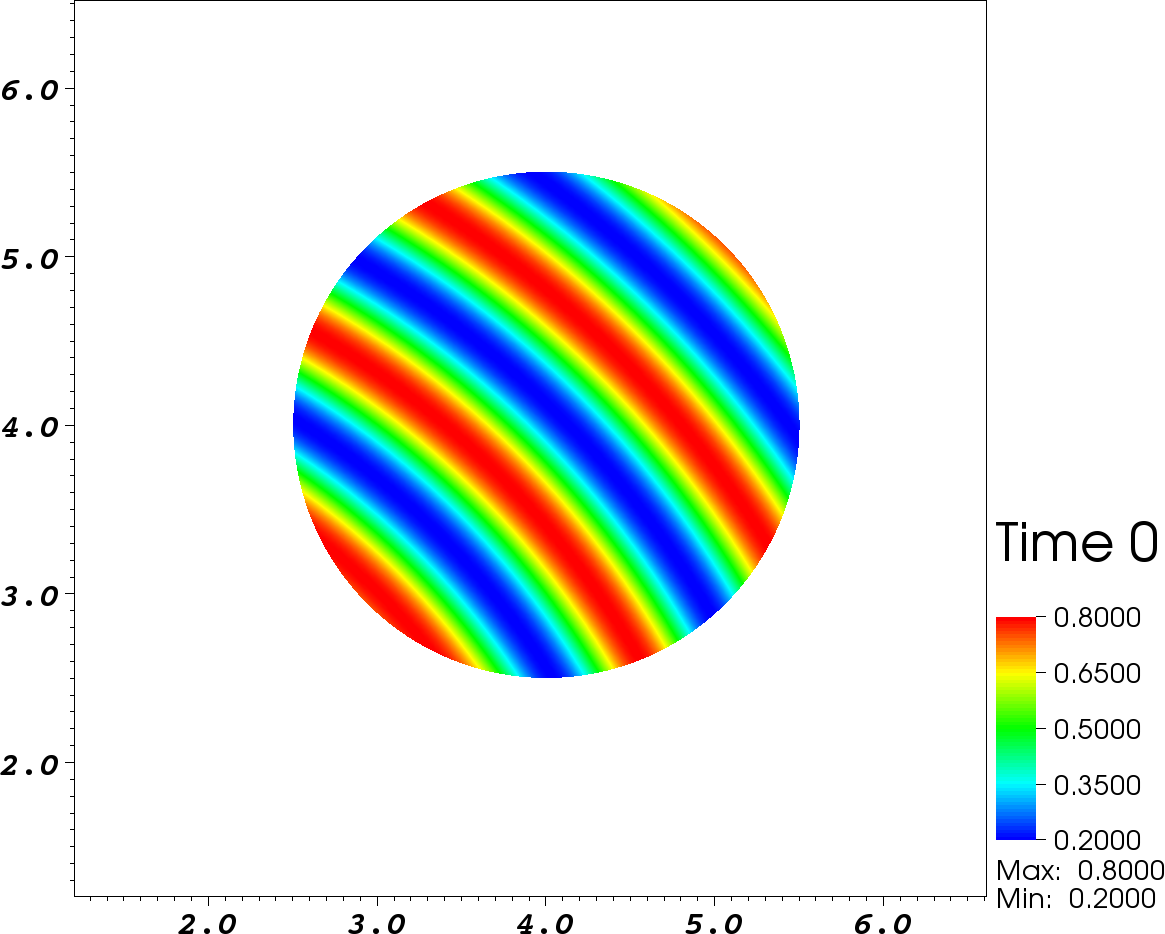}}
\subfloat[\label{subfig.}]{
	\includegraphics[width=.5\textwidth]{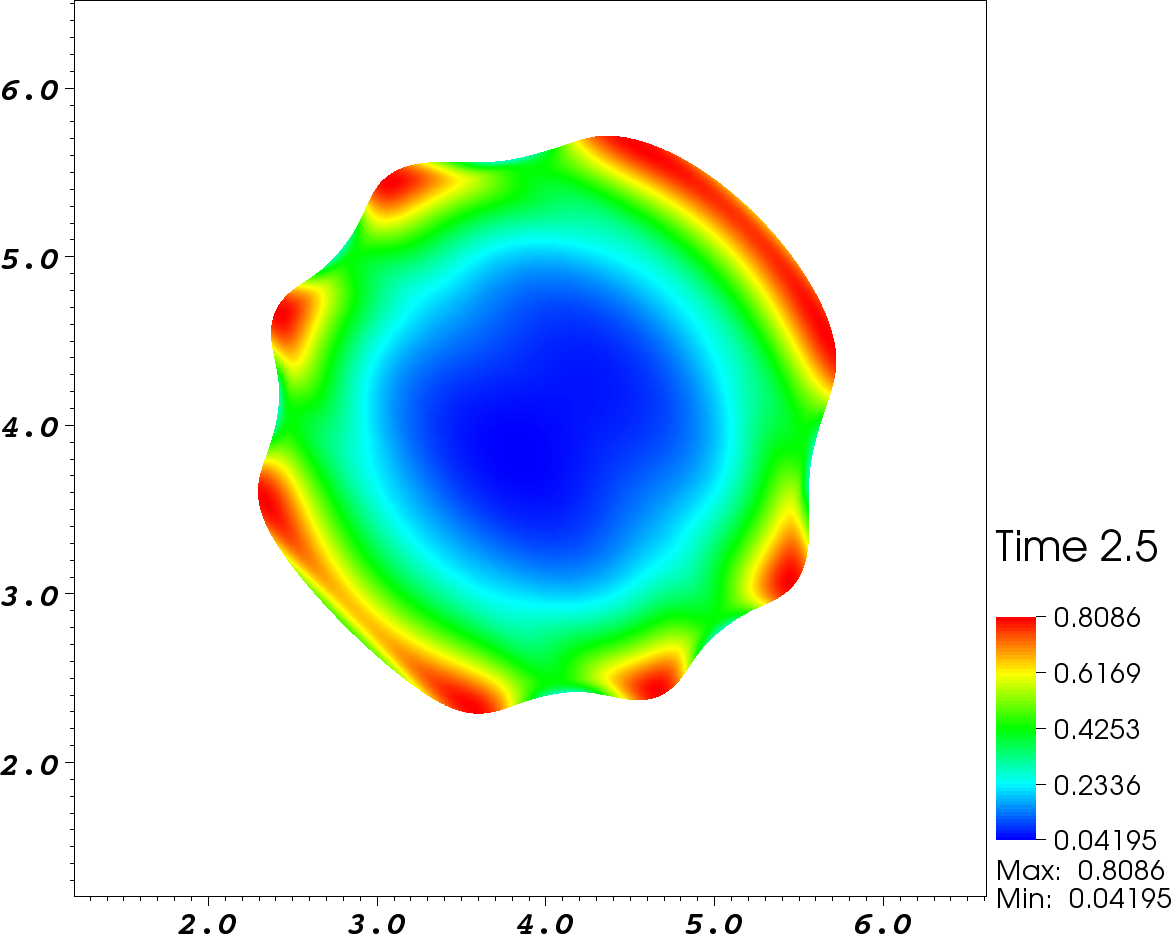}}\\
\subfloat[\label{subfig.}]{
	\includegraphics[width=.5\textwidth]{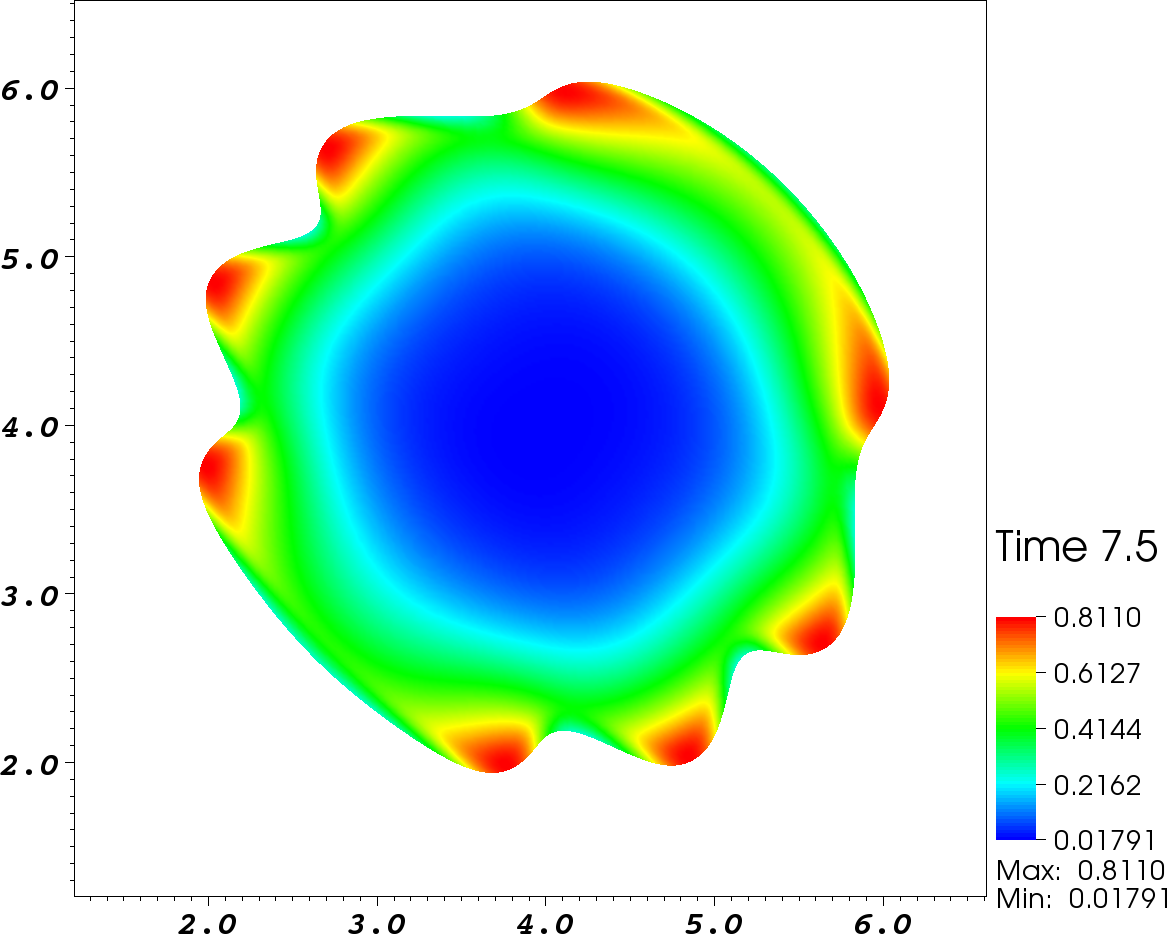}}
\subfloat[\label{subfig.}]{
	\includegraphics[width=.5\textwidth]{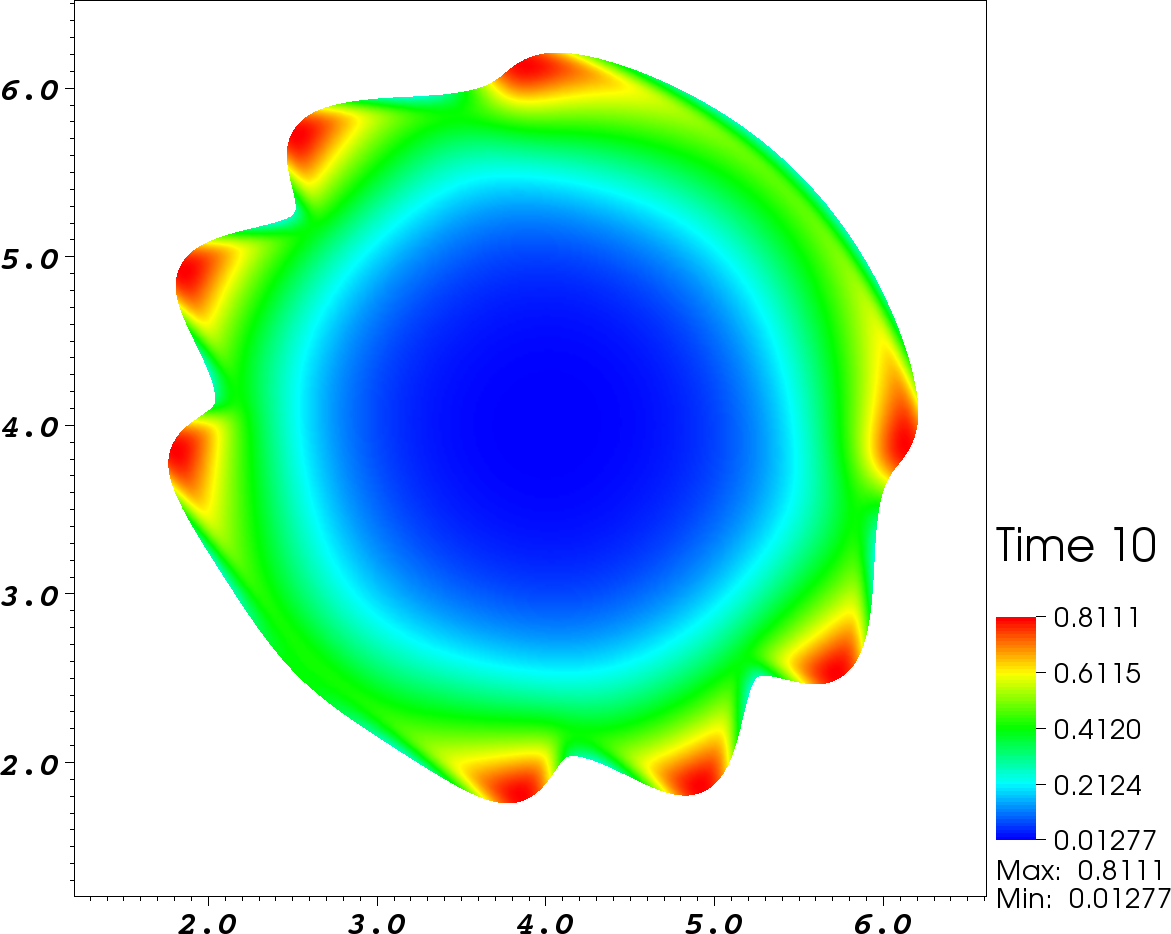}}
\end{center}
\caption{ECM distributions \dtr{at times: (a) $t=0$; (b) $t=2.5$; (c) $t=7.5$; and (d) $t=10$.}}
\label{distributions ECM}
\end{figure}

\section{Numerical results}
\label{num results}
In this section we show some numerical results obtained with the numerical method explained above. 
We have used the following initial conditions for cancer cells and ECM
\begin{align}
\label{initial conditions}
c_0(x) &= 
\left\{
\begin{array}{ll}
\displaystyle\frac{R - \|(x_1,x_2)-(4,4)\|_2}{R} & if\quad \|(x_1,x_2)-(4,4)\|_2 < R,\\
0 & else,
\end{array}
\right.\\
v_0(x) &=  0.3\, \sin\left(2 \pi \|(x_1,x_2)-(0,0)\|_2\right) + 0.5,
\end{align}
where $R$ is the initial radius of the cancer region and the point $(4,4)$ is the center of the computational domain, see Table \ref{tab:num setting} for the numerical parameters.

We first consider the configuration described in Table \ref{tab:conf1}. In particular, we consider a constant value for the proliferation coefficient
\begin{equation}
\label{const mu1}
\mu_1(v) = \mu_1^*,
\end{equation}
where $\mu_1^*$ is a constant shown in Table \ref{tab:conf1}.

\begin{figure}[!h]
\begin{center}
\subfloat[\label{subfig.}]{
	\includegraphics[width=.5\textwidth]{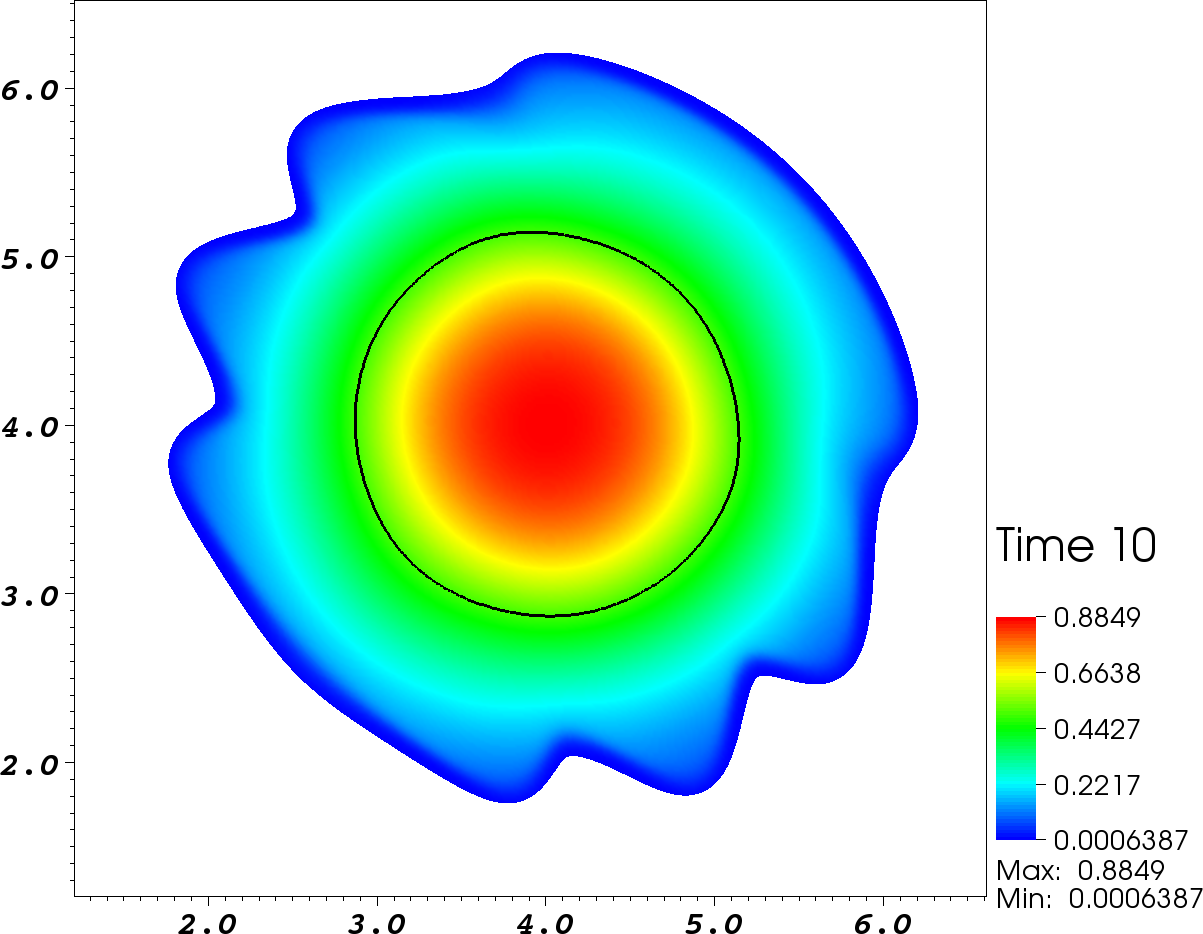}}
\subfloat[\label{subfig.}]{
	\includegraphics[width=.5\textwidth]{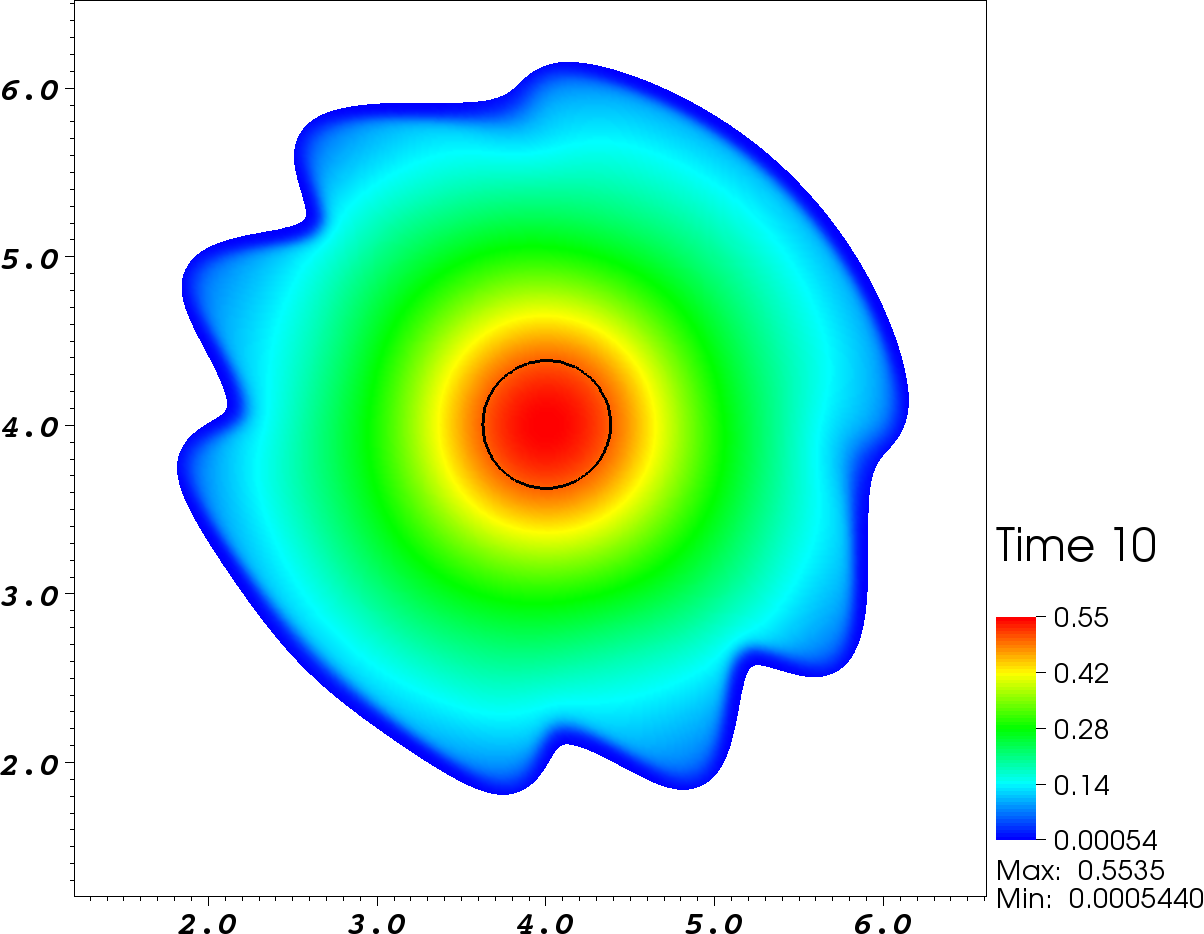}}
\end{center}
\caption{\dtr{Comparison of cancer distributions: (a) for a constant proliferation parameter $\mu_1 = \mu_1^*$; and (b) for ECM$-$dependent proliferation parameter $\mu_{1}(v)$ given in \eqref{mu_1}}. The black line is the contour of the level set $c\leq 0.5$.}
\label{distributions comparison}
\end{figure}

\tho{
It can be observed in Figure \ref{distributions cancer} that the cancer cells spread in the surrounding environment with a preferential path along the regions where the extracellular matrix has larger values. Due to fast degradation of the ECM in the central part of the cancer region, see Figure \ref{distributions ECM}, the transport of cancer cells following the gradient of ECM is limited to the vicinity of the boundary $\partial \Omega(t)$. The ECM shows some boundary layers at the interface with regions where the value of ECM outside of $\Omega(t)$ is larger. The layers are shown in red in the electronic version of the manuscript.
}

In the second \dtr{numerical} test, we \dtr{consider the ECM$-$dependent case for the proliferation parameter $\mu_{1}$ that was briefly introduced in Section \ref{section21_macro_model}. In this context, the proliferation parameter $\mu_1(v)$ is given here by 
\begin{equation}
\label{mu_1}
\mu_1(v) = 
\left\{
\begin{array}{ll}
\mu_1^*\, \exp\left(1+\displaystyle\frac{1}{(1-v)^2-1}\right) & if\quad v\in (0,1],\\
0 & if\quad v=0,
\end{array}
\right.
\end{equation}
which explores the proliferation conditions offered by the ECM, progressing smoothly from the worst conditions due to lack of ECM to the optimal conditions offered by abundant ECM density. }

\tho{
In Figure \ref{distributions comparison} we compare on the left side the distribution of cancer cells at time $t=10$ computed with constant $\mu_1 = \mu_1^*$ and on the right side the distribution computed with the function $\mu_1$ as in \eqref{mu_1}. The black line shown on both pictures is the isoline for the level $c(10, x)=0.5$. \dtr{While a similar morphology of the tumour boundary $\partial \Omega(10)$ is observed}, the effect of a decreasing $\mu_1$ with a decreasing $v$ results in \dtr{the spatial spread of cancer cells of high distribution levels (above $0.5$) being much more reduced in Figure \ref{distributions comparison}(b) than in Figure \ref{distributions comparison}(a), which corresponds to the case of constant $\mu_{1}$.} 
}

\section{Conclusions}
\label{conclusions}
We have presented a new formulation of a two-scale model to simulate cancer invasion. \dtr{This included a new derivation of the tumour boundary movement law by accounting on the MDE microdynamics contributions within a transport equation whose solution provides the level-set that indicates the new macro-scale cancer boundary.} 

 \dtr{At the core of the proposed numerical approach for }the proposed model \dtr{stands} the moving boundary method based on a combination of the level-set approach and cut-cells. The later are based on a continuous Lagrangian finite element formulation that on one side has high flexibility and accuracy properties and on the other side has revealed some instabilities issues that have been solved with a modification of the formulation. In particular, we have suppressed the cut-cells that were below a certain volume threshold. This numerical aspect is important and we are investigating the reason for such behaviour. Nevertheless, the neglection of the contribution of small cells can be interpreted as a quadrature error. By keeping the threshold for the suppression small, we keep this quadrature error small. 

\dtr{For the computational implementation, w}e have used an unfitted regular mesh with uniform cell diameters to avoid the problem of remeshing in case of large deformations.

A further important aspect for the numerical solution of the problem, that might be studied in a future work, is the singularly perturbed character of the macroscopic problem. 
Even if, in the considered configurations, we have not observed instability problems, it is of importance on its own to study possible stabilization techniques in combination with cut-cells.

We have shown that the presented framework is highly flexible to study possible variations of the model. In particular, since it is important to study the interplay between the two scales, the presented implementation allows high flexibility in defining the strength of the coupling via the definition of the velocity field.
In conclusion, we underline the potential of the presented method, that allows to go to three dimensional problems without changing the numerical formulation, \dtr{enabling this way a major development of the multiscale modelling framework introduced in \cite{trucu2013multiscale}.} The major \dtr{extension} needed for this development is the formulation of cut-cells in three dimensions.

\newpage
\appendix
\section{Parameter values}

\begin{table}[!h]
\begin{center}
\begin{tabular}{ccc}
  Final Time & T & 10\\
  \hline
  Initial radius of cancer distribution & $R$ & 1.5\\
  \hline
  Scale factor & $\varepsilon$ & 0.01\\
  \hline
  Diffusion cancer cells & $D_1$ & 0.0043\\
  \hline
  Convection & $\eta$ & 0.06\\
  \hline
  Proliferation & $\mu_1^*$ & 0.25\\
  \hline
  ECM remodelling & $\mu_2$ & 0.15\\
  \hline
  Degradation & $\alpha$ & 1.5\\
  \hline
  Diffusion \dtr{MDE} & $D_2$ & 0.001
\end{tabular}
\end{center}
\caption{Model parameters \dtr{in our numerical experiments}}
\label{tab:conf1}
\end{table}

\begin{table}[!h]
\begin{center}
\begin{tabular}{ccc}
  \tho{Splitting step} & $\Delta T$ & 0.1\\
  \hline
  Time step & $k$ & 0.1\\
  \hline
  Mesh size & $h$ & 0.015625\\
  \hline
  Stream-line stabilization & $\delta$ & 0.5\\
  \hline
  Radius of influence & $\rho$ & 0.1\\
  \hline
  Computational domain & $Y$ & $(0,8) \times (0,8)$ \\
  \hline
  \tho{Scaling factor for velocity $V(m)$} & $\cv$ & $5000$\\
  \hline
\end{tabular}
\end{center}
\caption{\dtr{Details} of \dtr{the} numerical setting}
\label{tab:num setting}
\end{table}

\section*{Acknowledgments}
TC was supported by the German Research Council (DFG) through project CA 633/2-1.
AVPB was funded by the Heidelberg Graduate School of Mathematical and Computational Methods for the Sciences (HGS MathComp), founded by DFG grant GSC 220 in the German Universities Excellence Initiative.

\end{document}

%% file: sampling.pdf_t
\begin{picture}(0,0)%
\includegraphics{sampling.pdf}%
\end{picture}%
\setlength{\unitlength}{4144sp}%
\begingroup\makeatletter\ifx\SetFigFont\undefined%
\gdef\SetFigFont#1#2#3#4#5{%
  \reset@font\fontsize{#1}{#2pt}%
  \fontfamily{#3}\fontseries{#4}\fontshape{#5}%
  \selectfont}%
\fi\endgroup%
\begin{picture}(4968,3596)(2736,-6351)
\put(3301,-3615){\makebox(0,0)[lb]{\smash{{\SetFigFont{20}{24.0}{\rmdefault}{\mddefault}{\updefault}{\color[rgb]{0,0,0}$\partial\Omega_h$}%
}}}}
\put(4289,-3155){\makebox(0,0)[lb]{\smash{{\SetFigFont{20}{24.0}{\rmdefault}{\mddefault}{\updefault}{\color[rgb]{0,0,0}$\varepsilon Y$}%
}}}}
\put(2751,-4491){\makebox(0,0)[lb]{\smash{{\SetFigFont{20}{24.0}{\rmdefault}{\mddefault}{\updefault}{\color[rgb]{0,0,0}$\varepsilon Y$}%
}}}}
\put(6678,-3022){\makebox(0,0)[lb]{\smash{{\SetFigFont{20}{24.0}{\rmdefault}{\mddefault}{\updefault}{\color[rgb]{0,0,0}$\varepsilon Y$}%
}}}}
\end{picture}%

%% file: cut-cells.pdf_t
\begin{picture}(0,0)%
\includegraphics{cut-cells.pdf}%
\end{picture}%
\setlength{\unitlength}{4144sp}%
\begingroup\makeatletter\ifx\SetFigFont\undefined%
\gdef\SetFigFont#1#2#3#4#5{%
  \reset@font\fontsize{#1}{#2pt}%
  \fontfamily{#3}\fontseries{#4}\fontshape{#5}%
  \selectfont}%
\fi\endgroup%
\begin{picture}(8402,4593)(2191,-5971)
\put(3967,-2029){\makebox(0,0)[lb]{\smash{{\SetFigFont{20}{24.0}{\rmdefault}{\mddefault}{\updefault}{\color[rgb]{0,0,0}\scalebox{2}{$\varphi_1$}}%
}}}}
\put(4520,-3689){\makebox(0,0)[lb]{\smash{{\SetFigFont{20}{24.0}{\rmdefault}{\mddefault}{\updefault}{\color[rgb]{0,0,0}\scalebox{2}{$\varphi_2$}}%
}}}}
\end{picture}%

%% file: extension.pdf_t
\begin{picture}(0,0)%
\includegraphics{extension.pdf}%
\end{picture}%
\setlength{\unitlength}{4144sp}%
\begingroup\makeatletter\ifx\SetFigFont\undefined%
\gdef\SetFigFont#1#2#3#4#5{%
  \reset@font\fontsize{#1}{#2pt}%
  \fontfamily{#3}\fontseries{#4}\fontshape{#5}%
  \selectfont}%
\fi\endgroup%
\begin{picture}(8492,4498)(2191,-5942)
\put(6174,-4987){\makebox(0,0)[lb]{\smash{{\SetFigFont{20}{24.0}{\rmdefault}{\mddefault}{\updefault}{\color[rgb]{0,0,0}\scalebox{1.8}{$t^n$}}%
}}}}
\put(8604,-4972){\makebox(0,0)[lb]{\smash{{\SetFigFont{20}{24.0}{\rmdefault}{\mddefault}{\updefault}{\color[rgb]{0,0,0}\scalebox{1.8}{$t^{n+1}$}}%
}}}}
\end{picture}%

%% file: extension2.pdf_t
\begin{picture}(0,0)%
\includegraphics{extension2.pdf}%
\end{picture}%
\setlength{\unitlength}{4144sp}%
\begingroup\makeatletter\ifx\SetFigFont\undefined%
\gdef\SetFigFont#1#2#3#4#5{%
  \reset@font\fontsize{#1}{#2pt}%
  \fontfamily{#3}\fontseries{#4}\fontshape{#5}%
  \selectfont}%
\fi\endgroup%
\begin{picture}(16422,5222)(2191,-5971)
\put(14091,-5029){\makebox(0,0)[lb]{\smash{{\SetFigFont{20}{24.0}{\rmdefault}{\mddefault}{\updefault}{\color[rgb]{0,0,0}\scalebox{1.8}{$t^{n+1}$}}%
}}}}
\put(6174,-4987){\makebox(0,0)[lb]{\smash{{\SetFigFont{20}{24.0}{\rmdefault}{\mddefault}{\updefault}{\color[rgb]{0,0,0}\scalebox{1.8}{$t^n$}}%
}}}}
\end{picture}%

%% file: twocuts.pdf_t
\begin{picture}(0,0)%
\includegraphics{twocuts.pdf}%
\end{picture}%
\setlength{\unitlength}{4144sp}%
\begingroup\makeatletter\ifx\SetFigFont\undefined%
\gdef\SetFigFont#1#2#3#4#5{%
  \reset@font\fontsize{#1}{#2pt}%
  \fontfamily{#3}\fontseries{#4}\fontshape{#5}%
  \selectfont}%
\fi\endgroup%
\begin{picture}(20685,12268)(-4784,-13196)
\put(15256,-1951){\makebox(0,0)[lb]{\smash{{\SetFigFont{20}{24.0}{\rmdefault}{\mddefault}{\updefault}{\color[rgb]{0,0,0}\scalebox{3}{$t_{n+1}$}}%
}}}}
\put(15886,-7576){\makebox(0,0)[lb]{\smash{{\SetFigFont{20}{24.0}{\rmdefault}{\mddefault}{\updefault}{\color[rgb]{0,0,0}\scalebox{3}{$t_n$}}%
}}}}
\put(-4769,-10546){\makebox(0,0)[lb]{\smash{{\SetFigFont{20}{24.0}{\rmdefault}{\mddefault}{\updefault}{\color[rgb]{0,0,0}\scalebox{3}{\bf cell cut twice}}%
}}}}
\end{picture}%